\documentclass[11pt]{article}

\usepackage{fullpage}
\usepackage{graphicx} 
\usepackage{amssymb,dsfont,color}
\usepackage{amsthm}
\usepackage{soul}
\usepackage{mathtools}

\usepackage{amsmath}
\usepackage{mathrsfs}
\usepackage{enumerate}
\usepackage{mathtools}
\usepackage{stmaryrd}

\usepackage{enumitem}

\def\s{\sigma}
\def\R{\mathbb{R}}
\def\e{\varepsilon}
\def\S{\cal{S}}

\def\I{\cal{I}}
\def\ds{\displaystyle}

\newtheorem{theorem}{Theorem}
\newtheorem{proposition}{Proposition}
\newtheorem{lemma}{Lemma}
\newtheorem{corollary}{Corollary}
 \newtheorem{DHassumption*}{$D$-Hurwitz assumption}
 
 \newtheorem{definition}{Definition}
 
 \newtheorem{remark}{Remark}

\begin{document}

\title{Stability criteria for singularly perturbed impulsive linear switched systems}

\author{Ihab~Haidar\thanks{Ihab Haidar is with Laboratoire d'\'Electronique et d'Automatique (L\'EA), Cergy, France, and with Laboratoire Jacques-Louis Lions (LJLL), Inria, Sorbonne Universit\'e, Paris, France, {\tt ihab.haidar@ensea.fr}.} \quad
Yacine~Chitour\thanks{Yacine Chitour is with Laboratoire des Signaux et Syst\`emes (L2S), Universit\'e Paris-Saclay, CNRS, CentraleSup\'elec, Gif-sur-Yvette, France,  {\tt yacine.chitour@l2s.centralesupelec.fr}.}\quad
Jamal~Daafouz\thanks{Jamal Daafouz is with CNRS, CRAN, Université de Loraine, Institut Universitaire de France (IUF), Vandoeuvre-lès-Nancy, France, {\tt Jamal.Daafouz@univ-lorraine.fr}.}\quad
Paolo~Mason\thanks{Paolo Mason \textcolor{blue}{is} with Laboratoire des Signaux et Syst\`emes (L2S), Universit\'e Paris-Saclay, CNRS, CentraleSup\'elec, Gif-sur-Yvette, France,  {\tt paolo.mason@l2s.centralesupelec.fr}.} \quad 
Mario~Sigalotti\thanks{Mario Sigalotti is with Laboratoire Jacques-Louis Lions (LJLL), Inria, Sorbonne Universit\'e, Universit\'e de Paris, CNRS, Paris, France, {\tt mario.sigalotti@inria.fr}.}
}

\maketitle

\begin{abstract}
We study a class of singularly perturbed impulsive linear switched systems exhibiting switching between slow and fast dynamics. To analyze their behavior, we construct auxiliary switched systems evolving in a single time scale. We prove that the stability or instability of these auxiliary systems directly determines that of the original system in the regime of small singular perturbation parameters.
\end{abstract}

%\begin{keyword}
%Switched systems, Hybrid systems, Singular perturbation, Converse Lyapunov, Exponential stability.
%\end{keyword}

 \section{Introduction}
 %\vspace{-0.3cm}
Consider the linear system evolving in $\R^d$
 %\vspace{-0.3cm}
\begin{equation*}\label{intro: eq1}
  \Sigma^{\e} :  \left\{
   \begin{array}{llll}
    {\cal D}^{\e}_k\dot X(t)&=&\Lambda_kX(t),\quad &t\in [t_k,t_{k+1}),\\
    X(t_{k+1})&=&R_k \ds\lim_{t\nearrow t_{k+1}}X(t), \quad  &k\ge 0,
    \end{array}\right.
     %\vspace{-0.3cm}
\end{equation*}
where $\Lambda_k,R_k$ take values in a compact subset of $d\times d$ real matrices and ${\cal D}^{\e}_k$ is a diagonal matrix with diagonal entries in $\{1,\e\}$, $\e$ being a small positive parameter. In this paper we deal with the problem of understanding the asymptotic behavior of this type of systems as $t$ goes to infinity in the regime where $\e$ is arbitrarily small.\\
System $\Sigma^{\e}$ represents a class of impulsive linear switched systems characterized by two time-scale dynamics: {\it fast variables}, whose velocities are modulated by $\frac1\e$, and the other ones called {\it slow variables}. The system is characterized by the dynamic interchange between slow and fast variables over time, governed by the switching signal $k\mapsto {\cal D}^{\e}_k$. While singularly perturbed hybrid systems with fixed slow-fast variables have been extensively studied in the literature (see, e.g., \cite{abdelgalil2023multi,9992853,simeonov1988stability,sanfelice2011singular,wang2012analysis,wang2025stabilization,Subotic2021,liberzon24,Tanwani24}), systems exhibiting switching slow/fast behaviors remain largely overlooked. Motivated by an industrial application in steel production (see~\cite{mallo09}), stability properties of the system $\Sigma^{\e}$ were first investigated in \cite{rejeb2018stability} in terms of LMI characterizations. As detailed in~\cite{mallo09}, steering control refers to strategies used to guide a metal strip through a finishing mill composed of $N$ stands, with the objective of keeping the strip aligned with the mill axis and avoiding excessive lateral deviations. During the last phase of a rolling process, called the tail-end-phase, the strip sequentially leaves the stands, inducing structural changes in the system dynamics. As a result, both the system dimension and the slow/fast decomposition vary across modes, which explains the mode-dependent singular perturbation structure. In particular, while the angles between the strip and the mill axis are usually fast variables, each time the strip leaves a stand, the angle at the first active stand becomes a slow state variable. In contrast to~\cite{rejeb2018stability}, existing results on switched singularly perturbed systems typically assume a fixed state dimension and a common time-scale separation across modes, and therefore do not apply in this context. Beyond hot strip mills, similar structures arise in systems with progressive loss of coupling or mode-dependent actuation, such as web handling or paper production processes. More generally, the proposed framework also applies to multi-agent systems with heterogeneous time scales, in which switching events may change the time-scale characteristics and role of individual agents.

The approach that we adopt here has been first explored in \cite{haidar2024necessary}, where some preliminary results to the present work were exhibited. Analyzing the stability of $\Sigma^{\e}$ is challenging as existing frameworks for singularly perturbed impulsive switched systems do not necessarily cover this class of systems. This work aims to address this gap by providing a comprehensive stability analysis for systems with switching slow/fast dynamics. It is important to emphasize that even in cases where the mappings $k\mapsto {\cal D}^{\e}_k$ and $k\mapsto R_k$ are constant, with
$R_k$ equal to the identity matrix, the classical singular perturbation theory \cite{KKO} cannot be applied in its standard form. In particular, the stability of the full system cannot be deduced directly from the stability of its individual components. To address this challenge, various stability criteria have been proposed in the literature (see, e.g., \cite{chitour2023upper,Hachemi2011,Malloci2009CDC,sanfelice2011singular}). For example, in \cite{chitour2023upper}, upper and lower bounds were derived for the maximal Lyapunov exponent of singularly perturbed linear switched systems as $\e$ tends to zero. In \cite{sanfelice2011singular}, stability was established under a dwell-time condition, which, importantly, does not explicitly depend on the time-scale parameter. Additionally, a recent study in \cite{TANG2024111462} explores the stabilization of switched affine singularly perturbed systems with state-dependent switching laws.\\
The purpose of this paper is twofold: first, to provide necessary or sufficient conditions ensuring a specific time-asymptotic behavior for $\Sigma^{\e}$ in the regime where $\e\sim0$, and second, to establish upper and lower bounds for the limit of the maximal Lyapunov exponent of $\Sigma^{\e}$ as $\e$ tends to $0$. Recall that the maximal Lyapunov exponent of a linear switched system represents the largest asymptotic exponential rate, as time tends to infinity, among all trajectories of the system. Stability conditions then emerge as special cases: specifically, a positive lower bound guarantees instability for all sufficiently small $\e$, while a negative upper bound ensures exponential stability for all $\e$ in a right-neighborhood of zero. This is provided after identifying some auxiliary discrete- and continuous-time single scale dynamics.\\ 
To carry out our analysis, we first rewrite system $\Sigma^{\e}$ in a new coordinate system that preserves the slow and fast nature of the variables over time. This is achieved through a mode-dependent variable reordering transformation, leading to a time-varying dimension for the slow and fast variables.
Starting from this new representation, we follow the classical Tikhonov approach to introduce auxiliary impulsive switched systems. In particular we introduce two continuous-time impulsive switched systems $\bar\Sigma$  and $\tilde\Sigma$ with reduced dimensions approximating the slow dynamics of $\Sigma^{\e}$. System $\bar\Sigma$ is obtained by neglecting the transient behavior during mode transitions while system $\tilde\Sigma$ is obtained by including the transient dynamics into the jump part of $\bar\Sigma$. Based on these two auxiliary systems and under suitable assumptions, 
we give bounds on the limit as $\e$ tends to $0$ of the maximal Lyapunov exponent of $\Sigma^{\e}$ as the following 
 %\vspace{-0.2cm}
\begin{equation}\label{sandwich}
\lambda(\bar\Sigma)\le \liminf_{\e\searrow0
}\lambda(\Sigma^{\e})\le 
\limsup_{\e\searrow0
}\lambda(\Sigma^{\e})\le \lambda(\tilde\Sigma). 
 %\vspace{-0.3cm}
\end{equation}
Observe that the left-hand side  inequality in \eqref{sandwich} yields a necessary condition for the stability of $\Sigma^{\e}$, in the sense that if $\bar\Sigma$ is exponentially 
unstable then there exists $\e_0>0$ such that for every $\e\in(0,\e_0)$ system $\Sigma^{\e}$ is exponentially unstable as well. On the other hand,  the right-hand side  inequality in \eqref{sandwich} yields a sufficient condition for the stability of $\Sigma^{\e}$, in the sense that if $\tilde\Sigma$ is exponentially stable then there exists $\e_0>0$ such that  for every $\e\in(0,\e_0)$ system $\Sigma^{\e}$ is exponentially stable as well. Under a dwell-time constraint, given that switching occurs slowly with respect to the time-scale $\frac1\e$, the transient phase is too short to aﬀect the dynamics of slow variables. Hence, in this case, systems $\bar\Sigma$ and $\tilde\Sigma$ 
have the same asymptotic behavior, 
leading to a complete characterisation of the limit as $\e$ tends to $0$ of the maximal Lyapunov exponent of $\Sigma^{\e}$. Another auxiliary single-scale dynamics denoted by $\hat\Sigma$ representing the transient behavior of $\Sigma^{\e}$ is also introduced. Based on $\hat\Sigma$, the limit as $\e$ tends to $0$ of the maximal Lyapunov exponent of $\Sigma^{\e}$ satisfies the inequality
%\vspace{-0.3cm}
\begin{equation}\label{1/e}
\lambda(\hat\Sigma)\le  \max\{0,\ds\liminf_{\e\searrow0}\e \lambda(\Sigma^{\e})\}, %\vspace{-0.3cm}
 \end{equation} 
giving a necessary condition for the stability of $\Sigma^{\e}$ in terms of $\hat\Sigma$. In fact, from~\eqref{1/e}, it follows that the exponential instability of $\hat\Sigma$ implies the exponential instability of $\Sigma^{\e}$ for every $\e>0$ sufficiently small. Moreover, 
 $\lambda(\Sigma^\e)$  is at least at order $\frac1\e$ as $\e$ tends to $0$.\\
The paper is organised as follows. In Section~\ref{sec:prob}, we reformulate system $\Sigma^\e$ within a suitable mathematical class and introduce the notion of stability for impulsive linear switched systems. We also state a stability theorem from~\cite{technical-spin-off} concerning the stability of impulsive linear switched systems, which serves as a central tool for the subsequent analysis. 
Section~\ref{sec:aux} introduces the auxiliary switched systems $\bar\Sigma_{\tau}$, $\hat\Sigma$, and $\tilde\Sigma$, and presents the main contributions through two 
theorems. The proofs of these theorems are detailed in Sections~\ref{sec: necessary} and~\ref{sec: sufficient}.
They rely on
a series of auxiliary results, provided in Section~\ref{sec:Lyapunov exponent}, enabling the reformulation of systems $\bar\Sigma_{\tau}$, $\hat\Sigma$, and $\tilde\Sigma$ within the impulsive switched system framework. Additional technical details are provided in the Appendix.
Section~\ref{sec:applications} addresses a particular class of $\Sigma^\e$ called the {\it complementary case}, and presents an illustrative example.

\subsection{Notation}\label{sec:notation}
 %\vspace{-0.3cm}
By $\R$ we denote the set of real numbers and by $\R_{\geq \tau}$ the set of real numbers greater than $\tau\geq 0$. We use $\mathbb{N}$ for the set of positive integers.
We use $M_{n,m}(\R)$ to denote the set of $n\times m$ real matrices and simply $M_n(\R)$ if $n=m$. The $n\times n$ identity matrix is denoted by $I_n$. 
By ${\rm GL}(n,\R)$ we denote the set of $n\times n$ invertible real matrices. 
For $Q\in M_{n,m}(\R)$ and $\ell \le n$, $\mathit{c}\le m$,
we denote by $(Q)_{\ell,\mathit{c}}$ the $\ell\times \mathit{c}$ matrix obtained by truncating $Q$ and keeping only its first $\ell$ lines and first $\mathit{c}$ columns. 
The spectral radius of a square matrix $M$ (i.e., the maximal modulus of its eigenvalues) is denoted by $\rho(M)$ and its spectral abscissa (i.e., the maximal real part of its eigenvalues) by $\alpha(M)$. \\
The Euclidean norm of a vector $x\in \R^n$ is denoted by $|x|$, while $\|\cdot\|$ denotes the induced norm on $M_n(\R)$, that is, $\|M\|=\max_{x\in \R^n\setminus\{0\}}\frac{|Mx|}{|x|}$ for $M\in M_n(\R)$. \\
Given $x:\R_{\ge 0}\to \R^n$ 
and $t>0$, we set $x(t^-):=\lim_{s\nearrow t}x(s)$ if such limit exists.\\
Given a set ${\cal Z}$, we denote by $\S_{{\cal Z}}$ the set of right-continuous piecewise-constant functions from $\R_{\geq 0}$ to ${\cal Z}$, that is, 
those functions $Z:\R_{\geq 0}\to {\cal Z}$ such that
there exists 
an increasing sequence $(t_k=t_k(Z))_{k\in \Theta^\star(Z)}$  of switching times in $(0,+\infty)$
which is 
locally finite (i.e., has no finite density point)  and for which $Z|_{[t_k,t_{k+1})}$ is constant 
for $k,k+1\in \Theta^\star(Z)$ (with $Z|_{[0,t_1)}$ and $Z|_{(\sup_{k\in \Theta^\star(Z)}t_k,+\infty)}$ also constant).  Here 
$\Theta^\star(Z)=\emptyset$, 
 $\Theta^\star(Z)=\{1,\dots,n\}$, 
or $\Theta^\star(Z)=\mathbb{N}$, 
depending on whether $Z$ has no, $n\in \mathbb{N}$, or infinitely many switchings, respectively.  Set $t_0=0$ and, when $\Theta^\star(Z)$ is finite with cardinality $n$, 
$t_{n+1}=+\infty$.  The value of $Z$ on $[t_k,t_{k+1})$ is denoted by $Z_k$.
Given $\tau\ge 0$, we denote by $\S_{{\cal Z},\tau}\subset \S_{{\cal Z},0}=\S_{\cal Z}$ the set of  piecewise-constant signals  with dwell time $\tau\ge 0$ (i.e., such that $t_{k+1}\ge t_k+\tau$ for $k\in \Theta(Z):=\{0\}\cup\Theta(Z)$). \\ 
Given a positive integer $d$, for every $\ell\in \{1,\dots, d\}$ and $\e \ge 0$ we use $E_{\ell}^{\e}$ to denote the $d\times d$ diagonal matrix with diagonal coefficients equal to $1$ over the $\ell$ first lines and $\e$ elsewhere. We denote by $E_{\ell^c}^{\e}$ the $d\times d$ diagonal matrix with diagonal coefficients equal to $\e$ over the $\ell$ first lines and $1$ elsewhere. Note that, for $\e>0$, we have $\e(E_{\ell}^{\e})^{-1}=E_{\ell^c}^{\e}$.

 %\vspace{-0.3cm}
\section{Problem formulation and main assumption}\label{sec:prob}
 %\vspace{-0.3cm}
\subsection{Singularly perturbed switched system}
 %\vspace{-0.3cm}
We begin this section by establishing a detailed reformulation of system $\Sigma^\e$. Let us fix an integer $d\ge 2$ and a compact subset ${\cal K}$ of $\{1,\dots, d-1\}\times {\rm GL}(d,\R)\times M_{d}(\R)\times M_{d}(\R)$.
We will use $\sigma$ to denote either an element of ${\cal K}$ or a 
signal in ${\cal S}_{\cal K}$, always specifying which case we are considering. The components of $\s$ will be denoted by 
$(\ell,P,\Lambda,R)$. In particular, if $\s$ is in ${\cal S}_{\cal K}$,
then $\ell$, $P$, $\Lambda$, and $R$ are themselves piecewise-constant signals 
and the  constant values taken by $\sigma$ are denoted by  $(\ell_k,P_k,\Lambda_k,R_k)$,  $k\ge 0$.
\\
For $\e>0$, $\tau\geq 0$ and $\sigma=(\ell,P,\Lambda,R)\in {\cal S}_{{\cal K},\tau}$, we introduce the system
\begin{equation*}\label{familly}
\Sigma_{\cal{K},\tau}^{\e} :
   \left\{
\begin{aligned}
    E^{\e}_{\ell_k}P_{k}\dot X(t)&=\Lambda_{k}X(t),  &\hspace{-0.7cm} t\in [t_k,t_{k+1}), \, k\in \Theta(\sigma), \\
    X(t_{k})&=R_{k-1} X(t_{k}^-),   &k\in \Theta^\star(\sigma),
\end{aligned}\right.
    \end{equation*}   
where $E^{\e}_{\ell_k}$ is a diagonal matrix with diagonal coefficients equal to $1$ over the $\ell_k$ first lines and $\e$ elsewhere.
The matrix $E^{\e}_{\ell_k}P_{k}$ identifies on each interval of time $[t_k,t_{k+1})$ the slow and fast variables of the system. The sets $\Theta(\s)$ and $\Theta^\star(\s)$, introduced in Section~\ref{sec:notation}, are used to parameterize
the switching instants of the signal $\sigma\in {\cal S}_{{\cal K},\tau}$. We denote by $\Phi_{\s}^{\e}(t,0)$ the flow at time $t$ of system $\Sigma_{\cal{K},\tau}^{\e}$ corresponding to the switching signal $\s\in {\cal S}_{{\cal K},\tau}$,
i.e., the matrix such that  $X_0\mapsto\Phi_{\s}^{\e}(t,0) X_0$ maps the initial condition $X(0)=X_0$  to the 
evolution at time $t$ of the
corresponding solution of $\Sigma_{\cal{K},\tau}^{\e}$. 
In analogy with the equality $\S_{{\cal K},0}=\S_{\cal K}$,
System $\Sigma_{{\cal K},0}^{\e}$ will be  denoted simply by $\Sigma_{\cal K}^{\e}$.
%\vspace{-0.3cm}
\begin{remark}
The case $\ell=d$, i.e., when all variables are slow, can be addressed by adding an extra fast variable $X_{d+1}$ in $\Sigma_{\cal{K},\tau}^{\e}$, %~\eqref{familly}, 
for example, defined by $\e \dot X_{d+1} = -X_{d+1}$. The stability analysis of this augmented system is equivalent to that of the original system, thus covering the case $\ell=d$.
\end{remark}
%\vspace{-0.3cm}
For a  fixed $\e>0$, $\Sigma_{{\cal K},\tau}^{\e}$ is a special case of the class of impulsive linear switched systems studied in \cite{technical-spin-off}. 
In next section we recall how such systems are defined and some crucial results 
about their
exponential stability.

\subsection{Impulsive linear switched systems}
 %\vspace{-0.3cm}
The definition of impulsive switched linear system and the main notions concerning its stability are recalled by the following definition. 
\begin{definition}\label{0-GES} 
Let $\tau\ge 0$, $d\in \mathbb{N}$, and ${\cal Z}$ be a bounded subset of  $M_d(\R)\times M_d(\R) $.
An \emph{impulsive switched linear system} is a switched system with state jumps of the form
\begin{equation*} 
    \Delta_{{\cal Z},\tau}: 
    \left\{\begin{aligned}
 \dot x(t)&=Z_1(t_k)x(t), & \hspace{-0.4cm} t\in [t_k,t_{k+1}), \, k\in \Theta(Z),\\
    x(t_{k})&= Z_2(t_{k-1})x( t_{k}^-), & k\in \Theta^{\star}(Z),\vspace{5pt}  
\end{aligned}\right.
\end{equation*} 
where $Z\in \S_{{\cal Z},\tau}$.
Denote by $\Phi_{Z}(t,0)$ the flow from time $0$ to time $t$ of $\Delta_{{\cal Z},\tau}$ corresponding to the switching signal $Z$. 
System $\Delta_{{\cal Z},\tau}$ is said to be
\begin{enumerate}
\item
\emph{exponentially  stable} (ES, for short) if there exist $c>0$ and $\delta>0$ such that 
%\begin{equation*}\label{eq-ES}
$\|\Phi_{Z}(t,0)\|\leq c e^{-\delta t}$
for every $t\geq 0$ and $Z\in {\cal S}_{{\cal Z},\tau};$
%\end{equation*}

\item
\emph{exponentially  unstable} (EU, for short) if there exist $c>0$, $\delta>0$, $Z\in {\cal S}_{{\cal Z},\tau}$, and $x_0\in\R^d\backslash\{0\}$ such that
$|\Phi_{Z}(t,0)x_0|\ge c e^{\delta t}|x_0|$ for every $t\ge 0$.
\end{enumerate}
The \emph{maximal Lyapunov exponent} of $\Delta_{{\cal Z},\tau}$ is defined as
\begin{equation*}%\label{eq:max-lyap-expo}
    \lambda(\Delta_{{\cal Z},\tau})=\limsup_{t\to+\infty}\sup_{Z\in {\cal S}_{{\cal Z},\tau}}\frac{\log(\|\Phi_{Z}(t,0)\|)}{t},
\end{equation*}
with the convention that $\log(0)=-\infty$. We define also the quantity $\mu(\Delta_{{\cal Z},\tau})$ given by
\begin{equation*}
\mu(\Delta_{{\cal Z},\tau})=\sup_{Z\in {\cal S}_{{\cal Z},\tau},\;k\in \Theta^\star(Z)}\frac{\log(\rho(\Phi_Z(t_k,0)))}{t_k}.
\end{equation*}
%\vspace{-0.5cm}
\end{definition}
Notice that for $\mu\in \R$, setting ${\cal Z}^{\mu}=\{(Z_1+\mu I_d,Z_2)\mid (Z_1,Z_2)\in {\cal Z}\}$, we have $\lambda(\Delta_{{\cal Z}^\mu,\tau})=\lambda(\Delta_{{\cal Z},\tau})+\mu$.

Let us introduce the notation $\Xi_{\cal Y}$ for a discrete-time switched system with set of modes ${\cal Y}\subset M_d(\R)$, that is,
\[\Xi_{{\cal Y}}:\quad x(k)=Y_k x(k-1),\qquad k\in \mathbb{N},\ Y\in {\cal Y}^\mathbb{N}.
\]
Recall that $\Xi_{{\cal Y}}$ is said to be \emph{bounded} if there exists a constant $C>0$ such that for every $k\in\mathbb{N}$ and every $Y_1,\dots,Y_k\in{\cal Y}$, $\|Y_k\cdots Y_1\|\leq C$. Otherwise, it is said to be  \emph{unbounded}. We will also say that $\Xi_{{\cal Y}}$ is \emph{exponentially unstable} (EU) if there exist $c>0,\delta>0,x_0\in\mathbb{R}^d\setminus\{0\}$, and a sequence of matrices $\{Y_k\}_{k\geq 0}$ in $\mathcal{Y}$ such that $\|Y_k\cdots Y_1 x_0\|\geq c e ^{\delta k}\|x_0\|$ for every $k\ge 1$.\

The next lemma provides an alternative characterization of the exponential stability of an impulsive linear switched system, formulated through its Lyapunov exponent.

\begin{lemma}[{\cite[Th. 1, Cor. 1 and Rem. 5]{technical-spin-off}}]\label{thm:supsup}
Let ${\cal Y}=\{Z_2\mid (Z_1,Z_2)\in {\cal Z}\}$.
Then $\lambda(\Delta_{{\cal Z},\tau})=+\infty$ if and only if $\tau=0$ and $\Xi_{{\cal Y}}$ is unbounded.
Moreover, if $\tau>0$ or system $\Xi_{{\cal Y}}$ is bounded,
then the following properties hold:
%\vspace{-0.5cm}
\begin{enumerate}
\item  \label{item1} $\lambda(\Delta_{{\cal Z},\tau}) = 
\displaystyle\max\left(\sup_{(Z_1,Z_2)\in {\cal Z}}\alpha(Z_1),\mu(\Delta_{{\cal Z},\tau})\right)$; 
\item \label{item2} $\Delta_{{\cal Z},\tau}$ is ES if and only if $\lambda(\Delta_{{\cal Z},\tau})<0$;
\item \label{item3} $\Delta_{{\cal Z},\tau}$ is EU if and only if $\lambda(\Delta_{{\cal Z},\tau})>0$.
\end{enumerate}
\end{lemma}

\subsection{Problem statement and first stability result}
 %\vspace{-0.3cm}
When $\e>0$ is fixed, $\Sigma^{\e}_{\cal K, \tau}$ is clearly an impulsive linear switched system.
Our goal is to characterize when $\Sigma^{\e}_{\cal K,\tau}$ is ES or EU for all values of $\e>0$ small enough, that is, according to Lemma~\ref{thm:supsup}, when $\lambda(\Sigma^{\e}_{\cal K,\tau})$ is negative or positive  for all values of $\e>0$ small enough.
A first trivial remark that can be done is that, since ${\cal S}_{{\cal K},\tau_1}\subset {\cal S}_{{\cal K},\tau_2}$ for $\tau_1\ge \tau_2$, then 
if $\Sigma^{\e}_{\cal K,\tau}$ is ES (respectively, EU) then $\Sigma^{\e}_{\cal K,\tilde\tau}$ is 
ES for every $\tilde\tau\in [\tau,+\infty)$
(respectively, EU for every $\tilde\tau\in [0,\tau]$). \\
In order to present some further remark on the exponential stability of $\Sigma^{\e}_{\cal K, \tau}$, let us 
introduce the following notation: given $\s=(\ell,P,\Lambda,R)\in{\cal K}$,  
we set
 %\vspace{-0.3cm}
\begin{equation}\label{ABCD}
\begin{pmatrix}
 A(\s) & B(\s)\\
 C(\s)& D(\s)
 \end{pmatrix}=
 \Lambda P^{-1},
 %\vspace{-0.3cm}
\end{equation}
where $A(\s)\in M_{\ell}(\R)$ and $B(\s),C(\s),D(\s)$ have the corresponding dimensions. 
When it is clear from the context, we simply write $A,B,C,D$ instead of $A(\s),B(\s),C(\s),D(\s)$.

 Let us introduce  
 ${\cal R}=\{R\mid (\ell,P,\Lambda,R)\in {\cal K}\}$
 and provide 
 a first  result, which uses the notation $E_{\ell^c}^0$ introduced in Section~\ref{sec:notation}.

\begin{proposition}\label{limit-criterium}
It holds that 
\begin{equation}\label{eq:limit-criterium}
    \liminf_{\e\searrow 0}\e\lambda(\Sigma^\e_{\cal K})\ge \sup_{(\ell,P,\Lambda,R)\in \cal K }\alpha(P^{-1}E_{\ell^c}^0\Lambda)\ge 0.
     %\vspace{-0.3cm}
\end{equation}
If, moreover, $\Xi_{{\cal R}}$ is bounded then 
 %\vspace{-0.3cm}
\[\lim_{\e\searrow 0}\e\lambda(\Sigma^\e_{\cal K})=\max(0,\lambda(\Delta_{{\cal Z},0})),
 %\vspace{-0.3cm}
 \]
where  ${\cal Z}=\{(P^{-1}E^{0}_{\ell^c}\Lambda,R)\mid (\ell,P,\Lambda,R)\in {\cal K}\}.$
\end{proposition}

\begin{proof} 
First notice that $\sup_{(\ell,P,\Lambda,R)\in \cal K }\alpha(P^{-1}E_{\ell^c}^0\Lambda)$ is nonnegative because it is larger than or equal to the real part of each eigenvalue of each matrix $P^{-1}E_{\ell^c}^0\Lambda=P^{-1}\left(\begin{smallmatrix}0&0\\C&D\end{smallmatrix}\right)P$, which is nonnegative because $\ell< d$.\\
Consider now,  for a given $\e\ge 0$, the impulsive linear switched system 
$\Delta_{{\cal Z}^\e,0}$
with ${\cal Z}^\e=\{(P^{-1}E^{\e}_{\ell^c}\Lambda,R)\mid (\ell,P,\Lambda,R)\in {\cal K}\}$. Then notice that, for every $\e>0$, 
the time rescaling $t\mapsto t/\e$ yields $\Phi_{Z}(t,0)=\Phi^{\e}_{\s}(\e t,0)$,
where $\s$ is an arbitrary signal in ${\cal S}_{\cal K}$ and $Z$ is 
the corresponding signal in ${\cal S}_{{\cal Z}^\e}$. This implies at once that
$\e\lambda(\Sigma^{\e}_{\cal K})=\lambda(\Delta_{{\cal Z}^\e,0})$.
Next, by Lemma~\ref{thm:supsup} we have that $\lambda(\Delta_{{\cal Z}^\e,0})\ge \sup_{(\ell,P,\Lambda,R)\in \cal K }\alpha(P^{-1}E^{\e}_{\ell^c}\Lambda)$.
It follows that $\e\lambda(\Sigma^{\e}_{\cal K
})\ge \sup_{(\ell,P,\Lambda,R)\in \cal K }\alpha(P^{-1}E_{\ell^c}^\e\Lambda)$. The proof of \eqref{eq:limit-criterium} is completed by letting $\e$ go to zero on both sides of the last inequality.\\
The last part of the statement comes from the fact that, if 
$\Xi_{{\cal R}}$ is bounded then  $\lambda(\Delta_{{\cal Z}^\e,0})<+\infty$ for every $\e\ge 0$ (Lemma~\ref{thm:supsup}). 
The convergence of $\lambda(\Delta_{{\cal Z}^\e,0})$ to $\lambda(\Delta_{{\cal Z}^{0},0})$ as $\e$ tends to $0$ is then a consequence of \cite[Proposition~6 and Remark~7]{technical-spin-off}. 
\end{proof}
Proposition~\ref{limit-criterium} immediately yields a sufficient condition for the exponential instability of 
$\Sigma^{\e}_{\cal K}$, namely that $\alpha(P^{-1}E_{\ell^c}^0\Lambda)>0$ for some $(\ell,P,\Lambda,R)\in {\cal K}$. 
Notice that for each $\s=(\ell,P,\Lambda,R)\in{\cal K}$ one has 
$P^{-1}E_{\ell^c}^0\Lambda=P^{-1}\left(\begin{smallmatrix}0&0\\C(\s)&D(\s)\end{smallmatrix}\right)P$.
This motivates the introduction of the following assumption. 

%\vspace{-0.1cm}
\begin{DHassumption*}%[$D$-Hurwitz assumption]
\label{Hurwitz assumption}
For each $\s\in{\cal K}$, the matrix $D(\s)$ defined in~\eqref{ABCD} is Hurwitz.
\end{DHassumption*}

%\vspace{-0.3cm}

\section{Auxiliary switched systems and statement of the main results}\label{sec:aux}
 %\vspace{-0.3cm}
The stability of $\Sigma_{\cal{K},\tau}^{\e}$ will be studied by comparing it with that of single-scale auxiliary systems, which are introduced in this section. 

\subsection{Block diagonalization}
%\vspace{-0.3cm}

Following a classical approach (see e.g. \cite{kokotovic1975riccati}), for $\s=(\ell,P,\Lambda,R)\in {\cal K}$ we introduce the transformation matrix $T^{\e}=T^{\e}(\s)$ given by
%\vspace{-0.3cm}
\begin{equation*}\label{Transformation}
T^{\e}=\left(\begin{smallmatrix}
I_{\ell}& 0\\
D^{-1}C+\e Q^{\e} & I_{d-\ell}
\end{smallmatrix}\right)P,
%\vspace{-0.3cm} 
\end{equation*}
and the upper triangular matrix $\Gamma^{\e}=\Gamma^{\e}(\s)$ given by
%\vspace{-0.3cm}
\begin{equation*}
\Gamma^{\e}=
\left(\begin{smallmatrix}
    A-BD^{-1}C-\e BQ^\e&B\\
    0&\frac{D}{\e}+  
    (D^{-1}C+\e Q^{\e})B
\end{smallmatrix}\right),
%\vspace{-0.3cm}
\end{equation*}
 where $Q^{\e}=Q^{\e}(\s)$ is chosen in such a way 
that
%\begin{equation*}
$\frac{1}{\e}T^\e P^{-1}E^\e_{\ell^c}\Lambda (T^\e)^{-1}=\Gamma^{\e}$,  
%\end{equation*} 
and $\|Q^\e(\s)\|$ is upper bounded uniformly with respect to 
$\s\in{\cal K}$  and $\e$ small enough.\\
Notice that the coordinate transformation just introduced makes sense only if the 
matrix $D$ is invertible. The proof of the existence of $Q^\e$ can be found in \cite{kokotovic1975riccati}. \\
Let us stress that 
the expression for $T^\e$ makes sense also for $\e=0$, and we will  write simply $T(\s)$ for $T^0(\s)$. Note that the matrices $T(\sigma)$ belong to a compact subset of invertible matrices.\\
The transformation above allows one to introduce 
the variables $x(t)$ and $z(t)$ of dimensions $\ell(t)$ and $d-\ell(t)$, respectively, such that
%\vspace{-0.3cm}
\begin{equation}\label{variable reordering}
\begin{pmatrix}
x(t)\\
z(t)
\end{pmatrix}
= T^{\e}_{k}X(t), \quad  \forall\,t\in [t_k,t_{k+1}), \, k\in \Theta(\s),
%\vspace{-0.3cm}
\end{equation}
and system~$\Sigma_{\cal{K},\tau}^{\e}$  
can be equivalently represented in terms of the triangular matrices $\Gamma^\e=\Gamma^\e(\s)$, for $\sigma\in {\cal S}_{{\cal K},\tau}$, 
    as
\begin{equation}\label{xzeqdot}
\left\{
\begin{aligned}
\begin{pmatrix}
\dot x(t)\\
\dot z(t)
\end{pmatrix}
&=
\Gamma^\e_k
\begin{pmatrix}
x(t)\\
z(t)
\end{pmatrix}, & 
\hspace{-0.3cm}
\begin{array}{l}
t\in [t_k,t_{k+1}),\\ k\in \Theta(\sigma) 
\end{array}\\
\begin{pmatrix}
x(t_{k})\\
z(t_{k})
\end{pmatrix}
&= T^{\e}_{k}R_{k-1}(T^{\e}_{k-1})^{-1}
\begin{pmatrix}
x(t_{k}^-)\\
z(t_{k}^-)
\end{pmatrix},& k\in \Theta^\star(\sigma).
\end{aligned}\right. 
%\vspace{-0.3cm}
\end{equation} 
Observe from~\eqref{xzeqdot} that even if the starting model~$\Sigma_{\cal{K},\tau}^{\e}$  
does not include jumps in its dynamics, i.e., if
$R= I_{d}$ for every $(\ell,P,\Lambda,R)\in\mathcal{K}$, the change of variables~\eqref{variable reordering} leads anyway to a singularly perturbed switched system with jumps.
 
%\vspace{-0.3cm}
\subsection{Slow dynamics by Tikhonov's approach}
 %\vspace{-0.3cm}
The Tikhonov decomposition of a singularly perturbed system consists in analyzing the limit behavior of the slow dynamics by setting $\e=0$ and replacing in the equation of the slow dynamics the limit  
value of the fast variable.
This can be done when the fast dynamics has a stable equilibrium (as long as the switching signal stays constant), that is, when $D(\s)$ is Hurwitz, for $\s\in{\cal K}$. 
Assuming that the $D$-Hurwitz assumption holds and applying this approach to \eqref{xzeqdot} leads to the formulation of the  
system 
%\vspace{-0.3cm}
\begin{equation*}\label{barsigma}
\bar\Sigma_{\tau}:  
\left\{\begin{aligned}
\dot{\bar x}(t) &= M_k\bar x(t), & t\in [t_k,t_{k+1}), \, k\in\Theta(\sigma) \vspace{5pt}\\
\bar x(t_{k}) 
&= J(k)\bar x(t_{k}^-), & k \in\Theta^\star(\s), 
\end{aligned}\right. 
%\vspace{-0.3cm}
\end{equation*} 
where $\s \in{\cal S}_{{\cal K},\tau}$, $M_k=A_k-B_kD_k^{-1}C_k$
and $J(k)=\left(T_{k}R_{k-1}T_{k-1}^{-1}\right)_{\ell_{k},\ell_{k-1}}$.\\
In what follows, we write $\bar \Sigma$ for $\bar\Sigma_0$. We also introduce the  subset of $M_d(\R)$  given by 
%\vspace{-0.3cm}
\begin{equation*}\label{barR}
\bar {\cal R}=\left\{RT^{-1}
\left(\begin{smallmatrix}
I_{\ell}&0\\
0&0
\end{smallmatrix}\right)
T\mid (\ell,P,\Lambda,R)\in {\cal K}\right\},
%\vspace{-0.3cm}
\end{equation*}
which is related to the jumps of system $\bar\Sigma_{\tau}$. \\
When dwell-time is active ($\tau>0$) the reduced system $\bar\Sigma_{\tau}$ allows to establish  a necessary and a sufficient condition for the stability of  $\Sigma^{\e}_{\cal K,\tau}$. More precisely, we will prove that 
if system $\bar\Sigma_{\tau}$ is ES (respectively, EU) then $\Sigma^{\e}_{\cal K,\tau}$ is ES (respectively, EU) for every $\e>0$ small enough (cf. Theorems~\ref{main theorem-necessary} and \ref{main theorem-sufficient}).

\subsection{Transient dynamics}
 %\vspace{-0.3cm}
If there is no  dwell-time constraint (i.e., $\tau=0$), the transient dynamics governed by the fast dynamics must be considered. To capture such transient dynamics, a rescaling of time is needed and new variables are introduced: $s=t/\e$, $\hat x(s)=x(\e s)$ and $\hat z(s)=z(\e s)$. After rewriting the dynamics~\eqref{xzeqdot} in terms of this new scale, the limit problem at $\e=0$ is given, for $\s=(\ell,P,\Lambda,R)\in{\cal S}_{\cal K}$, by 
\begin{equation*}\label{fast}
\hat\Sigma:  \left\{
\begin{aligned}
\begin{pmatrix}
\dot{\hat x}(s)\\
\dot{\hat z}(s)
\end{pmatrix}
&=
\begin{pmatrix}
0 & 0\\
0 & D_k
\end{pmatrix}  
\begin{pmatrix}
\hat x(s)\\
\hat z(s)
\end{pmatrix},
& \hspace{-0.3cm}
\begin{array}{l}
s \in [s_k, s_{k+1}) \\
k \in \Theta(\sigma)
\end{array}
\vspace{5pt}\\
\begin{pmatrix}
\hat x(s_{k})\\
\hat z(s_{k})
\end{pmatrix} 
&=T_{k}R_{k-1}T^{-1}_{k-1} 
 \begin{pmatrix}
\hat x(s_{k}^-)\\
\hat z(s_{k}^-)
\end{pmatrix} , &k \in\Theta^\star(\sigma).
\end{aligned}\right.
\vspace{-0.3cm}
\end{equation*} 
System $\hat\Sigma$ allows one to establish a necessary condition for the stability of $\Sigma^{\e}_{\cal K}$, in the sense that its instability implies the instability of $\Sigma^{\e}_{\cal K}$ (Theorem~\ref{main theorem-necessary}). \\
Finally, we introduce an additional auxiliary system, 
obtained from $\bar\Sigma$ by allowing more complex jumps, which take into account 
the transient dynamics
 described by $\hat\Sigma$. Consider
\begin{equation*}\label{tildesigma}
\tilde\Sigma:  
\left\{\begin{aligned}
\dot{\tilde x}(t) &= M_k\tilde x(t), & t\in [t_k,t_{k+1}), \, k\in\Theta(\sigma) \vspace{5pt}\\
\tilde x(t_{k})
&= \tilde J(k)
\tilde x(t_{k}^-)
& k \in \Theta^\star(\s),
\end{aligned}\right.
\end{equation*} 
for $\s\in{\cal S}_{\cal K}$, where $\tilde J(k)=\left(T_{k} F_{k-1} R_{k-1}T_{k-1}^{-1}\right)_{\ell_{k},\ell_{k-1}}$
and $F_{k-1}$ is any element of  $\hat{\cal F}$, where
\begin{align*}
\hat {\cal F} = \{I_d\} \cup \bigg\{ 
&\prod_{i=1}^{n} R_i T_i^{-1}
\left(\begin{smallmatrix}
I_{\ell_{i}} & 0 \\
0 & e^{s_i D_{i}}
\end{smallmatrix}\right)
T_i \ \bigg| \ n\in\mathbb{N}, s_i > 0, 
(\ell_i, P_i, \Lambda_i, R_i)\in {\cal K}
\text{ for } i = 1,\dots,n \bigg\}.
\end{align*}
Intuitively speaking, $\tilde\Sigma$ takes into account at once the two cases in which the difference between subsequent switching times is much larger or comparable to the parameter $\varepsilon$.\\
Although $\bar\Sigma_{\tau}$, $\hat\Sigma$ and $\tilde\Sigma$ are not formally impulsive linear switched system in the sense of  Definition~\ref{0-GES} (their jump dynamics depend at each time $t_k$ on the value of $\s$ both on $[t_{k},t_{k+1})$ and $[t_{k-1},t_{k})$) and that their state dimensions may vary with time (as in the case of $\bar\Sigma_{\tau}$ and $\tilde\Sigma$), their stability properties can be defined in analogy with Definition~\ref{0-GES}.
In particular, for $\Sigma=\bar\Sigma_{\tau}$ ($\hat\Sigma$, $\tilde\Sigma$, respectively), we denote by $\Phi^{\Sigma}_\sigma(t,0)$ the flow from time $0$ to time $t\ge 0$ of $\Sigma$ associated with a signal $\s\in {\cal S}={\cal S}_{\cal K,\tau}$ (${\cal S}_{\cal K}$, respectively) and introduce
%\vspace{-0.3cm}
\begin{equation}\label{lambdasigma}
\lambda(\Sigma)=\limsup_{t\to+\infty}\sup_{\s\in {\cal S}}\frac{\log(\|\Phi^{\Sigma}_{\s}(t,0)\|)}{t}.
%\vspace{-0.3cm}
\end{equation}
We also introduce for $\hat\Sigma$ the Lyapunov-like exponent
%\vspace{-0.3cm}
\begin{equation*} \tilde\lambda(\hat\Sigma)=\limsup_{s\to+\infty}
\sup_{\s\in {\cal S}_{\cal K},\;k\in \Theta^\star(\s),\;s=s_k}\frac{\log(\|\Phi^{\hat\Sigma}_\s(s,0)\|)}{s},
%\vspace{-0.3cm}
\end{equation*}
obtained by considering the evolution only at switching times.\\ 
The relation between the exponents defined above and the corresponding stability properties  will be discussed in detail in Section~\ref{sec:Lyapunov exponent}.

\subsection{Main results}
 %\vspace{-0.3cm}
Our main results are summarized in the following two theorems. 
The first one contains, in particular, 
conditions under which 
$\Sigma^{\e}_{\cal K,\tau}$ is EU for every $\e$ small enough.

\begin{theorem}\label{main theorem-necessary}
Assume that the $D$-Hurwitz assumption holds true. The following statements hold:
\begin{enumerate}

\item\label{itemii} 
For every\, $\tau>0$, we have 
%\begin{eqnarray}
$\lambda(\bar\Sigma_\tau)\le \liminf_{\e\searrow0}\lambda(\Sigma^{\e}_{\cal K,\tau})$. 
%\label{eq:inequality2}
 %\end{eqnarray} 
If, moreover, $\bar\Sigma_\tau$ is EU then $\Sigma^{\e}_{\cal K,\tau}$ is EU for every $\e>0$ small enough.

\item\label{itemiii} If $\tau=0$ and both $\Xi_{{\cal R}}$ and $\Xi_{\bar {\cal R}}$ are bounded then the inequality stated in item~\ref{itemii} holds true and if, moreover, $\bar\Sigma$ is EU then $\Sigma^{\e}_{\cal K}$ is EU for every $\e>0$ small enough.

\item\label{itemv} 
 If\, $\Xi_{\cal R}$ is bounded,  then
%\begin{eqnarray}
 $\lambda(\hat\Sigma)\le  \max\{0,\ds\liminf_{\e\searrow0}\e \lambda(\Sigma^{\e}_{\cal K})\}$.
 %\label{eq:inequality1}
%\end{eqnarray} 
In particular, if $\hat\Sigma$ is EU then,  for every $\e>0$ small enough, $\Sigma^{\e}_{\cal K}$ is EU and
 $\lambda(\Sigma^\e_{\cal K})$  is 
at least at order $1/\e$ as $\e$ tends to $0$.
\end{enumerate}

\end{theorem}

The second theorem collects results containing sufficient conditions for the ES of 
$\Sigma^{\e}_{\cal K,\tau}$ for $\e$ small enough. 

\begin{theorem}\label{main theorem-sufficient}
Assume that the $D$-Hurwitz assumption holds true. 
The following statements hold:
\begin{enumerate}
\item\label{itemi} 
For every $\tau>0$, we have 
%\begin{eqnarray}
$\lambda(\bar\Sigma_{\tau})\ge \lim_{\e\searrow0}\lambda(\Sigma^{\e}_{\cal K,\tau})$.
 %\label{eq:inequality0}
 %\end{eqnarray} 
 In particular, if $\bar\Sigma_\tau$ is ES then $\Sigma^{\e}_{\cal K,\tau}$ is ES for every $\e>0$ small enough.

\item\label{itemiv} 
Assume that $\tilde\lambda(\hat\Sigma)<0$. 
Then, 
%\begin{eqnarray}
$\limsup_{\e\searrow0}\lambda(\Sigma^{\e}_{\cal K})\le \lambda(\tilde\Sigma)$. 
% \label{eq:inequality3}
 %\end{eqnarray} 
In particular, if $\tilde\Sigma$ is ES then $\Sigma^{\e}_{\cal K}$ is ES for every $\e>0$ small enough. 

\end{enumerate}

\end{theorem}

As a direct consequence of Theorems~\ref{main theorem-necessary} and \ref{main theorem-sufficient}, we obtain the following corollary. 
\begin{corollary}\label{<>=}
Assume that the $D$-Hurwitz assumption holds true and that $\tau>0$. 
Then  
%\begin{equation*}
$\lambda(\bar\Sigma_{\tau})=\lim_{\e\searrow0}\lambda(\Sigma^{\e}_{\cal K,\tau})$.
 %\label{eq:equality000}
 %\end{equation*} 
\end{corollary}
In the simplified case of switched singular perturbations with constant $\ell$ and $P,R\equiv I_{d}$, the corollary takes the following form, which completes the results obtained in \cite{chitour2023upper}.

\begin{corollary}\label{automatica}
Let $\tau>0$ and ${\mathcal M}$ be a compact subset of $M_d(\R)$. 
Consider the singularly perturbed linear switched system 
\begin{equation*}
\Upsilon^{\e}_{\tau}:\left\{
\begin{aligned}
\dot x(t)&=A_k x(t)+B_k y(t),\\
\e\dot y(t)&=C_k
x(t)+D_k
y(t),
\end{aligned}
\right.  t\in[t_k,t_{k+1}),\ k\in\Theta(M),
\end{equation*}
for $M=(\begin{smallmatrix}A&B\\C&D\end{smallmatrix})\in {\cal S}_{{\cal M},\tau}$. Suppose that $D$ is Hurwitz for every $(\begin{smallmatrix}A&B\\C&D\end{smallmatrix})\in {\cal M}$. 
Consider the reduced order system 
\begin{equation*} 
\bar\Upsilon_{\tau}: \dot {\bar x}(t)= M_k\bar x(t), \quad t\in[t_k,t_{k+1}),\ k\in\Theta(M), 
\end{equation*}
where $M_k=A_k-B_kD_k^{-1}C_k$, for $M=(\begin{smallmatrix}A&B\\C&D\end{smallmatrix})\in {\cal S}_{{\cal M},\tau}$. Then  
%\begin{equation*}
 $\lambda(\bar\Upsilon_{\tau})=\lim_{\e\searrow0
 }\lambda(\Upsilon^{\e}_{\tau}).$
%\end{equation*} 
In particular, if $\bar\Upsilon_{\tau}$ is ES (respectively, EU) then $\Upsilon^{\e}_\tau$ is ES (respectively, EU) for every $\e>0$ small enough. 
\end{corollary}
The proofs of Theorems~\ref{main theorem-necessary} and 
\ref{main theorem-sufficient} are provided in Sections~\ref{subsec:itemii}, \ref{subsec:itemiii},  \ref{subsec:itemv}, \ref{subsec:itemi}, and \ref{subsec:itemiv}. In order to obtain 
these proofs, we introduce a series of results that allow us to reformulate the auxiliary switched systems $\bar\Sigma_{\tau}$, $\hat\Sigma$, and $\tilde\Sigma$ in the framework of impulsive switched systems, as defined in Definition~\ref{0-GES}. This reformulation is essential for properly characterizing the Lyapunov exponents of the auxiliary switched systems, which exhibit time-varying dimensions and  ``multi-mode-dependent" jumping parts. These preliminary results are developed in Section~\ref{sec:Lyapunov exponent}.
Moreover, additional technical preliminaries are required. In particular, Section~\ref{sec: necessary} presents an approximation result that links the flow of the singularly perturbed system to the auxiliary systems introduced. Additionally, Section~\ref{sec: sufficient} includes converse-type theorems for impulsive switched systems, originally established in~\cite{technical-spin-off} and adapted to our context, which play a key role in the proofs of Theorems~\ref{main theorem-necessary} and \ref{main theorem-sufficient}.

\section{Lyapunov exponents of singularly perturbed and auxiliary systems}\label{sec:Lyapunov exponent}

Given $\mu\in\R$, 
$\s=(\ell,P,\Lambda,R)\in {\cal K}$, and $\e>0$, we introduce
\begin{align}\label{Gamma}
\hspace{-0.3cm}M^\mu(\s)=M(\s)+\mu I_{\ell} \; \text{ and }  \; \Gamma^{\e,\mu}(\s)=\Gamma^{\e}(\s)+\mu I_{d}.
\end{align}
By $\Sigma^{\e, \mu}_{\cal K, \tau}$ we denote the $\mu$-shifted system associated with~\eqref{xzeqdot} and corresponding to $\Gamma^{\e,\mu}$. 
Notice that 
$\Gamma^{\e,0}$ coincides with $\Gamma^{\e}$
and
$\Sigma^{\e, 0}_{\cal K, \tau}$, up to the choice of coordinates, with $\Sigma^{\e}_{\cal K, \tau}$.  Notice also that $\lim_{\e\searrow 0}\alpha(\Gamma^{\e,\mu})=\alpha(M^\mu)$ for every $\s\in{\cal K}$.\\
The goal of this section is to relate the 
maximal Lyapunov exponent of $\Sigma_{\cal K,\tau}^{\e,\mu}$ and the 
exponential growth rates of the auxiliary systems associated with $\Sigma_{\cal K,\tau}^{\e}$, 
as introduced in \eqref{lambdasigma},
with those of suitably identified 
\emph{weighted discrete-time switched systems, as introduced in \cite{technical-spin-off}.
Let us recall some notation from that reference.} 
Given $n\in \mathbb{N}$ and a subset ${\cal N}$ of $M_n(\R)\times \R_{\ge 0}$, 
we denote by $\Omega_{\mathcal{N}}$  the set of all sequences $\omega=((N_j,\tau_j))_{j\in\mathbb{N}}$ in ${\cal N}$ such that $\sum_{j\in \mathbb{N}}\tau_j=+\infty$.  
For every $\omega=((N_j,\tau_j))_{j\in\mathbb{N}}\in \Omega_{\mathcal{N}}$ and $k\in \mathbb{N}$, we set
\[\omega_k=((N_j,\tau_j))_{j=1}^k,\, 
|\omega_k|=\tau_1+\dots+\tau_k,\,
\Pi_{\omega_k}=N_k\cdots N_1.\]
Given $\mu\in\R$, $\tau\ge 0$, and $\e>0$, 
we define 
\begin{align*}
{\cal N}^{\e,\mu}_{\tau}&=\left\{
\left(R{(T^{\e})}^{-1}
e^{t\Gamma^{\e,\mu}}
T^{\e},t\right)
 \mid  \s\in {\cal K}, t\geq \tau \right\},
 \end{align*}
 which is a subset of $M_d(\R)\times\R_{\geq 0}$.
We  denote ${\cal N}^{\e,0}_{\tau}$ simply
by ${\cal N}^{\e}_{\tau}$, 
${\cal N}^{\e,\mu}_{0}$ simply
by ${\cal N}^{\e,\mu}$, and ${\cal N}^{\e,0}_{0}$ simply
by ${\cal N}^{\e}$.

\begin{lemma}\label{lambda-sigmaeps}
If $\Xi_{{\cal R}}$ is unbounded then $\lambda(\Sigma^{\e}_{\cal K})=+\infty$ for every $\e>0$. 
If $\tau>0$ or 
system $\Xi_{{\cal R}}$ is bounded, then, for every $\mu\in \R$,  
\begin{align*}
\lambda(\Sigma^{\e,\mu}_{\cal K, \tau})&=
\max\left(\sup_{\s\in \cal K}\alpha(\Gamma^{\e,\mu}),\sup_{\substack{\omega\in \Omega_{{\cal N}^{\e,\mu}_{\tau}}\\ k\in \mathbb{N}}}\frac{\log(\rho(\Pi_{\omega_k}))}{|\omega_k|}\right)\\
&<+\infty.
\end{align*}
\end{lemma}
%\vspace{-0.3cm}
\begin{proof}
First notice that $\mu$ can be taken equal to zero, since $\lambda(\Sigma^{\e,\mu}_{\cal K, \tau})=\lambda(\Sigma^{\e}_{\cal K, \tau})+\mu$ and the terms in the right-hand side of \eqref{lambda-sigmaeps}
scale analogously.

Observe that 
$\Sigma^{\e}_{\cal K, \tau}$ can be equivalently written as 
\begin{equation*}
\left\{
\begin{aligned}
\dot X(t)
&=
(T^{\e}_k)^{-1}\Gamma_k^{\e}T^{\e}_k
X(t),& t\in[t_k,t_{k+1}), \, k\in \Theta(\sigma), \vspace{5pt}\\
X(t_{k})
&= R_{k-1}X(t_{k}^-),& k\in\Theta^\star(\sigma),
\end{aligned}\right.
\end{equation*}
for $\sigma=(\ell,P,\Lambda,R)\in {\cal S}_{{\cal K},\tau}$, that is, as the impulsive linear switched system $\Delta_{{\cal Z},\tau}$, with 
\[{\cal Z}=\{((T^{\e}(\s))^{-1}\Gamma^{\e}(\s)T^{\e}(\s),R(\s))\mid \s\in{\cal K}\}.\] 
Applying Lemma~\ref{thm:supsup} we have
that if $\Xi_{{\cal R}}$ is unbounded then $\lambda(\Sigma^{\e}_{\cal K})=+\infty$ for every $\e>0$, while if $\Xi_{{\cal R}}$ is bounded or $\tau>0$ 
then
\begin{align*}
&+\infty>\lambda(\Sigma^{\e}_{\cal K,\tau})=\\
&\max\left(\sup_{\s\in \cal K}\alpha((T^{\e})^{-1}\Gamma^{\e}T^{\e}),
\sup_{\s\in {\cal S}_{{\cal K},\tau},\;k\in \mathbb{N}}\frac{\log(\rho(\Phi^{\e}_\s(t_k,0)))}{t_k}\right)\\
&=\max\left(\sup_{\s\in \cal K}\alpha(\Gamma^{\e}),\sup_{\omega\in \Omega_{{\cal N}^{\e}_{\tau}}, k\in \mathbb{N}}\frac{\log(\rho(\Pi_{\omega_k}))}{|\omega_k|}\right),
\end{align*} 
concluding the proof. 
\end{proof}

We introduce the  subset of $M_d(\R)\times\R_{\geq 0}$ given by 
\begin{align*}\label{Nhat}
{\hat {\cal N}}&=\left\{\left(RT^{-1}\left(\begin{smallmatrix} 
I_{\ell}&0\\
0 & e^{s D(\s)} 
\end{smallmatrix}\right)T,s\right) \mid  \s\in {\cal K}, s> 0\right\}.
\end{align*}

\begin{lemma}\label{lambda-hatsigma}
If\, $\Xi_{{\cal R}}$ is unbounded then $\lambda(\hat\Sigma)=+\infty$. On the other hand, if $\Xi_{{\cal R}}$ is bounded, then 
\begin{equation}\label{eq:lambda-hatsigma}
\lambda(\hat\Sigma)=\max\left(0,\sup_{\omega\in \Omega_{\hat{\cal N}}, k\in \mathbb{N}}\frac{\log(\rho(\Pi_{\omega_k}))}{|\omega_k|}\right)<+\infty,
\end{equation}
and $\hat\Sigma$ is EU if and only if  $\lambda(\hat\Sigma)>0$.
\end{lemma}

\begin{proof}
Consider 
\begin{equation*} \label{hatDelta}
\hat\Delta:\left\{\begin{aligned}
&\dot X(s)
=
T^{-1}_{k}
\left(\begin{smallmatrix}
0 & 0\\
0 & D_k
\end{smallmatrix}\right) 
T_{k}
X(s),
 &
 \begin{array}{l}
  s\in [s_k,s_{k+1}),\\   k \in \Theta(\sigma), 
 \end{array} 
  \vspace{5pt}\\
&X(s_{k})
=R_{k-1} X(s_{k}^-) , & k \in\Theta^\star(\sigma),
\end{aligned}\right.
\end{equation*} 
for $\s=(\ell,P,\Lambda,R)\in{\cal S}_{\cal K}$. First notice that $\hat \Delta$ is an impulsive linear switched system in the sense of Definition~\ref{0-GES}. Moreover we have that the trajectories of $\hat \Sigma$ and $\hat \Delta$ only differ by a mode-dependent change of coordinates belonging to a compact set of invertible matrices, given as follows: for $t\geq 0$ there exists $\sigma\in {\cal{K}}$ so that 
$
X(t)=T(\sigma)^{-1}\left(\begin{smallmatrix}\hat x(t)\\ \hat z(t)\end{smallmatrix}\right)$.
As a consequence, we have that $\lambda(\hat\Sigma)=\lambda(\hat\Delta)$ 
and $\hat\Sigma$ is EU if and only if the same is true for $\hat \Delta$. \\
Applying Lemma~\ref{thm:supsup}, we have  that
if $\Xi_{{\cal R}}$ is unbounded then $\lambda(\hat\Delta)=+\infty$, while 
if $\Xi_{{\cal R}}$ is bounded
then
\begin{align*}
&+\infty>\lambda(\hat\Delta)=\\
&\max\left(\sup_{\s\in \cal K}\alpha\left(T^{-1}
\left(\begin{smallmatrix}
0&0\\
0&D(\s)
\end{smallmatrix}\right)T\right),\sup_{\omega\in \Omega_{\hat{\cal N}}, k\in \mathbb{N}}\frac{\log(\rho(\Pi_{\omega_k}))}{|\omega_k|}
\right)\\
&=\max\left(0,\sup_{\omega\in \Omega_{\hat{\cal N}}, k\in \mathbb{N}}\frac{\log(\rho(\Pi_{\omega_k}))}{|\omega_k|}\right)
\end{align*} 
and $\hat\Delta$ is EU if and only if  $\lambda(\hat\Delta)>0$,
concluding the proof. 
\end{proof}

\begin{remark}\label{rmk:compati}
Lemma~\ref{lambda-hatsigma} implies, in particular, that  $\lambda(\hat\Sigma)\ge 0$.
It is still possible 
that  $\tilde\lambda(\hat\Sigma)=\tilde\lambda(\hat\Delta)<0$. In that case, one has restrictions on the jumps, namely, $R(\s)x=0$ for every $x\in \R^d$ such that 
$\left(
\begin{smallmatrix} 
0 & 0\\ 
C(\s) & D(\s)
\end{smallmatrix}
\right)
P(\s)x=0 $, for every $\s\in {\cal K}$, as proved in \cite[Proposition~5]{technical-spin-off}. 
\end{remark}

For $\mu\in\R$ and $\tau\ge 0$, we introduce the  subset of $M_d(\R)\times\R_{\geq 0}$ given by 
\begin{align*}
\bar{\cal N}^{\mu}_{\tau}&=
\left\{
\left(RT^{-1}\left(\begin{smallmatrix} 
e^{t M^{\mu}} 
&0\\
0 & 0
\end{smallmatrix}\right)T,t\right)
 \mid  \s\in {\cal K},\; t\geq \tau \right\}.
\end{align*}
We  denote $\bar{\cal N}^{0}_{\tau}$ 
simply by $\bar{\cal N}_{\tau}$. Moreover, consider the system $\bar\Sigma_\tau^\mu$  built as $\bar \Sigma_\tau$ where we replace the matrix $M_k$ by the matrix $M^{\mu}_k=M_k+\mu I_{\ell_{k}}$.

\begin{lemma}\label{lambda-Sigmabar}
If\, $\Xi_{\bar{\cal R}}$ is unbounded then $\lambda(\bar\Sigma)=+\infty$. On the other hand,
if $\tau>0$ or 
system $\Xi_{\bar{\cal R}}$ is bounded, then
\begin{equation}\label{eq:lambda-Sigmabar}
\lambda(\bar\Sigma_\tau^\mu)=
\max\left(\sup_{\s\in \cal K}\alpha(M^\mu),\sup_{\substack{\omega\in \Omega_{\bar{\cal N}^\mu_{\tau}}\\ k\in \mathbb{N}}}\frac{\log(\rho(\Pi_{\omega_k}))}{|\omega_k|}\right)<+\infty,
\end{equation}
and $\bar\Sigma_\tau^\mu$ is ES (respectively, EU) if and only if $\lambda(\bar\Sigma_\tau^\mu)<0$ (respectively, $\lambda(\bar\Sigma_\tau^\mu)>0$).
\end{lemma}

\begin{proof}
As in the proof of Lemma~\ref{lambda-sigmaeps}, we assume without loss of generality that $\mu=0$. \\
For each $\s\in {\cal K}$, let us define $\bar R=\bar R(\s)$ as  $\bar R=RT^{-1}
\left(\begin{smallmatrix}
I_{\ell} & 0\\
0 & 0
\end{smallmatrix}\right)
T$.
Observe that, for every $\delta\in\R$,  $\bar{\cal N}_{\tau}$ can be equivalently written as 
\begin{align*} 
\bar{\cal N}_{\tau}=
\left\{
\left(\bar R
e^{t 
T^{-1}
\left(\begin{smallmatrix} 
M&0\\
0 & \delta I_{d-\ell}
\end{smallmatrix}\right)T}
,t\right)
 \mid  \s\in {\cal K}, t\geq \tau \right\}.
 \end{align*} 
Consider the associated linear impulsive system 
\begin{equation*}\label{barsigma-extended}
\bar\Delta_{\tau}^{\delta}:
\left\{\begin{aligned}
\dot X(t) &= T_k^{-1}
%\begin{pmatrix} 
\left(\begin{smallmatrix}
M_k&0\\
0 & \delta I_{d-\ell_k}
%\end{pmatrix}
\end{smallmatrix} \right)
T_k X(t), 
&\hspace{-0.3cm}
\begin{array}{l} 
t\in [t_k,t_{k+1}),\\
k\in \Theta(\s), 
\end{array}\\
X(t_{k})
&= \bar R_{k-1}
X(t_{k}^-), & k \in \Theta^\star(\s),
\end{aligned}\right.
\end{equation*} 
for $\s\in{\cal S}_{{\cal K},\tau}$, and let us denote by $\Phi_{\s}^{\bar\Delta}(t,0)$ the corresponding flow at time $t$ associated with a signal $\s\in {\cal S}_{\cal K}$. \\
If $\Xi_{\bar{\cal R}}$ is unbounded then, by
Lemma~\ref{thm:supsup}, $\lambda(\bar\Delta_{0}^{\delta})=+\infty$.
Actually, 
applying \cite[Lemma~3]{technical-spin-off} 
one can find sequences $(\omega^n)_{n\in \mathbb{N}}$ in $\Omega_{\bar{\cal N}_0}$ and $(k_n)_{n\in \mathbb{N}}$ in $\mathbb{N}$
such that $\lim_{n\to \infty} |\omega^n_{k_n}| =+\infty$ and
\begin{equation}\label{eq:limtilds}
\limsup_{n\to+\infty}\dfrac{\log(\|\Pi_{\omega^n_{k_n}}\|)}{|\omega^n_{k_n}|}
=+\infty.
\end{equation}
Now, observe that for every $\tau\ge0$, $\omega\in\Omega_{ {\bar{\cal N}_\tau}}$, and $k\in \mathbb{N}$
we have 
 \begin{equation}\label{cont-discret1}
\Pi_{\omega_k}
=R_kT_{k}^{-1}
\begin{pmatrix}
\Phi^{\bar\Sigma}_{\s_\omega}(|\omega_k|^-,0) &0\\
0&0
\end{pmatrix}
T_0
\end{equation}
where we recall that
$\Phi^{\bar\Sigma}_{\s_\omega}(|\omega_k|^-,0)$ denotes the limit as $s\nearrow |\omega_k|$ of the flow of system $\bar\Sigma_{\tau}$ from time $0$ to time $s$, associated with the signal $\s_\omega$ that corresponds to $\omega$.
We can then deduce from \eqref{eq:limtilds} and \eqref{cont-discret1} that \[\limsup_{n\to+\infty}\dfrac{\log(\|\Phi^{\bar\Sigma}_{\s_{\omega^n}}(|\omega^n_{k_n}|^-,0)\|)}{|\omega^n_{k_n}|} =+\infty,\]
yielding
$\lambda(\bar\Sigma)=+\infty$.\\
If either $\Xi_{\bar{\cal R}}$ is bounded or $\tau>0$, we deduce from Lemma~\ref{thm:supsup} that 
\begin{align*}
\lambda(\bar\Delta^{\delta}_{\tau}) =\max\left(\delta, \sup_{\s\in \cal K}\alpha(M),\sup_{\omega\in \Omega_{\bar{\cal N}_{\tau}},k\in\mathbb{N}}\frac{\log(\rho(\Pi_{\omega_k}))}{|\omega_k|}\right).
\end{align*}
Now, fix $\delta<\sup_{\s\in \cal K}\alpha(M)$ and let us prove that $\lambda(\bar\Sigma_\tau)=\lambda(\bar\Delta^{\delta}_{\tau})$.  
According to \cite[Proposition~4]{technical-spin-off} 
we can characterize $\lambda(\bar\Delta^{\delta}_{\tau}) $ also as
\[\lambda(\bar\Delta^{\delta}_{\tau}) =\max\left(\sup_{\s\in\cal K}\alpha(M),\sup_{\omega\in \Omega_{\bar{\cal N}_{\tau}}}\limsup_{k\to\infty}\frac{\log(\|\Pi_{\omega_k}\|)}{|\omega_k|}\right).\]
By consequence,  
using \eqref{cont-discret1}, we have
\begin{align*}
\lambda(\bar\Delta^{\delta}_{\tau}) &\le
\max\left(\sup_{\s\in\cal K}\alpha(M),\sup_{\omega\in \Omega_{\bar{\cal N}_{\tau}}}\limsup_{k\to\infty}\frac{\log(\|\Phi_{\s_\omega}^{\bar\Sigma}(|\omega_k|^-,0)\|)}{|\omega_k|}
\right)\\
&\leq \lambda(\bar\Sigma_{\tau}),
\end{align*}
where the last inequality follows from 
the definition of $\lambda(\bar\Sigma_{\tau})$ 
(considering a constant signal to deduce that $\alpha(M%(\s)
)\le \lambda(\bar\Sigma_{\tau})$ for every $\s\in {\cal K}$). \\
On the other hand, for $\s\in {\cal S}_{{\cal K},\tau}$,  $k\in \mathbb{N}$, and $t\in [t_k,t_{k+1})$ we have 
\begin{align*}
\Phi_{\s}^{\bar\Delta}(t,0)=T_{k}^{-1}
\begin{pmatrix}
\Phi_{\s}^{\bar\Sigma}(t,0)&0\\
\star & 0  
\end{pmatrix}
T_0
\end{align*}
from which we get the  inequality 
\begin{align}\label{eq:DeltaboundsSigma}
\|\Phi^{\bar\Sigma}_{\s}(t,0)\|\le 
\left\|\left(\begin{smallmatrix}
\Phi^{\bar\Sigma}_{\s}(t,0)&0\\
\star & 0  
\end{smallmatrix}\right)\right\|\le 
C\|\Phi^{\bar\Delta}_{\s}(t,0)\|,\quad \forall\;t\ge t_1,
\end{align}
for some $C>0$ depending only on ${\cal K}$. From the definition of Lyapunov exponent of $\bar\Sigma_{\tau}$ and $\bar\Delta_{\tau}^{\delta}$ it follows that 
\begin{align*}
\lambda(\bar\Sigma_{\tau})&=\limsup_{t\to+\infty}\sup_{\s\in {\cal S}_{{\cal K},\tau}}\frac{\log(\|\Phi^{\bar\Sigma}_{\s}(t,0)\|)}{t}\\
&\le 
\limsup_{t\to+\infty}\sup_{\s\in {\cal S}_{{\cal K},\tau}}\frac{\log(\|\Phi^{\bar\Delta}_{\s}(t,0)\|)}{t}= \lambda(\bar\Delta_{\tau}^{\delta}),
\end{align*}
concluding the proof that $\lambda(\bar\Sigma_\tau)=\lambda(\bar\Delta^{\delta}_{\tau})$. \\
Notice that  
$\lambda(\bar\Sigma_\tau^\mu)=\lambda(\bar\Delta^{\delta,\mu}_{\tau})$, where  $\bar\Delta^{\delta,\mu}_{\tau}$ is the natural shifted version of $\bar\Delta^{\delta}_{\tau}$. \\
Let us conclude the proof by showing  that, under the assumption that $\Xi_{\bar{\cal R}}$ is bounded, 
 $\bar\Sigma_\tau^\mu$ is ES (respectively, EU) if and only if $\lambda(\bar\Sigma_\tau^\mu)<0$ (respectively, $\lambda(\bar\Sigma_\tau^\mu)>0$).
One implication being trivial by definition of $\lambda(\bar\Sigma_\tau^\mu)$, 
let us assume that $\lambda(\bar\Sigma_\tau^\mu)<0$ (respectively, $\lambda(\bar\Sigma_\tau^\mu)>0$) and prove that $\bar\Sigma_\tau^\mu$ is ES (respectively, EU).
On the one hand, if $\lambda(\bar\Sigma_\tau^\mu)<0$ then, since 
$\lambda(\bar\Sigma_\tau^\mu)=\lambda(\bar\Delta^{\delta,\mu}_{\tau})$ and thanks to Lemma~\ref{thm:supsup}, 
 $\bar\Delta^{\delta,\mu}_{\tau}$ is ES. The exponential stability of $\bar\Sigma_\tau^\mu$ follows then from \eqref{eq:DeltaboundsSigma} and the fact that, according to \eqref{eq:lambda-Sigmabar}, $\alpha(M^\mu(\s))<0$ for every $\s\in {\cal K}$. 
On the other hand, if $\lambda(\bar\Sigma_\tau^\mu)>0$ then \eqref{eq:lambda-Sigmabar}
immediately identifies a constant or periodic signal yielding the exponential instability of 
$\bar\Sigma_\tau^\mu$.
 \end{proof}
 
 \begin{remark}
As a consequence of Lemma~\ref{lambda-Sigmabar}, if $\Xi_{\bar {\cal R}}$ is bounded, then $\lim_{\tau\searrow 0}\lambda(\bar\Sigma_\tau)=\lambda(\bar\Sigma)$.
Indeed, first notice that if $\Xi_{\bar {\cal R}}$ is bounded then $\lambda(\bar\Sigma_\tau)$ is characterized by \eqref{eq:lambda-Sigmabar} for every $\tau\ge 0$. 
Notice also that $\lambda(\bar\Sigma)\ge \lambda(\bar\Sigma_\tau)$ for every $\tau>0$. 
 On the other hand, if $\lambda(\bar\Sigma)=\sup_{\s\in \cal K}\alpha(M)$ then 
 $\lambda(\bar\Sigma)\le \lambda(\bar\Sigma_\tau)$ for every $\tau>0$. 
 We are left to prove that $\liminf_{\tau\searrow 0}\lambda(\bar\Sigma_\tau)\ge \lambda(\bar\Sigma)$ when for every $\delta>0$ 
 there exist $\omega\in \Omega_{\bar{\cal N}}$ and $k\in \mathbb{N}$ such that
 $\lambda(\bar\Sigma)\le \frac{\log(\rho(\Pi_{\omega_k}))}{|\omega_k|}+\delta$. 
 In that case,  
 for every $\tau>0$ small enough $\omega_k$ can be completed to a $k$-periodic sequence in  $\Omega_{\bar{\cal N}_\tau}$, so that $\lambda(\bar\Sigma_\tau)\ge \frac{\log(\rho(\Pi_{\omega_k}))}{|\omega_k|}\ge \lambda(\bar\Sigma)-\delta$, completing the proof of the claim.
 \end{remark}

\begin{remark}\label{lemma:checkbar}
Another consequence of Lemma~\ref{lambda-Sigmabar} is that, if $\Xi_{\bar{\cal R}}$ is EU, then there exists $\tau>0$ such that $\bar\Sigma_{\tau}$ is EU.
Indeed, 
if $\Xi_{\bar{\cal R}}$ is EU then
there exist $\bar R(0),\dots,\bar R(L-1)\in \bar{\mathcal{R}}$ such that $\rho(\bar R(L-1)\cdots \bar R(0))>1$ (see, e.g., \cite{jungers2009joint}). Let $\tau>0$ and $\omega\in \Omega_{\bar {\cal N}_\tau}$ be the $L$-periodic sequence given by $\Pi_{\omega_{L}}=\bar R(L-1)e^{\tau \bar M_{L-1}}\cdots \bar R(0)e^{\tau \bar M_0}$, where $\bar M_k=e^{T_{k}\left(\begin{smallmatrix}
M_k & 0\\
0& 0
\end{smallmatrix}\right)T_{k}}$ and $\bar R(k)=R_k T_k^{-1}\left(\begin{smallmatrix}
I_{\ell_k} & 0\\
0& 0
\end{smallmatrix}\right)T_{k}$, for $k=0,\dots, L-1$. 
By the continuity of the spectral radius, it follows that 
$\rho(\Pi_{\omega_{L}})>1$ for $\tau$ small enough.
Observe that $\sup_{\omega\in \Omega_{\tilde{\cal N}_{\tau}}, k\in \mathbb{N}}\frac{\log(\rho(\Pi_{\omega_k}))}{|\omega_k|}\ge \frac{\log(\rho(\Pi_{\omega_L}))}{L\tau}>0$, the conclusion follows from Lemma~\ref{lambda-Sigmabar}. 
\end{remark}
We introduce the  subset of $M_d(\R)\times\R_{\geq 0}$ given by 
\begin{align*}
\tilde{\cal N}^\mu&=
\left\{\left(F
  R T^{-1 }
 \left(\begin{smallmatrix}
 %\left(\begin{smallmatrix}
e^{tM^\mu}&0\\
0&0
 \end{smallmatrix}\right)
%\end{smallmatrix}\right)
 T,t\right)
 \mid  \s\in {\cal K},F\in \hat{\cal F},t\ge 0 \right\}
\end{align*}
and the  subset of $M_d(\R)$  given by 
\begin{equation*}
 \tilde {\cal R}=\left\{F RT^{-1}
\left(\begin{smallmatrix}
I_{\ell}&0\\
0&0
\end{smallmatrix}\right)
T \mid 
\s\in {\cal K},\;
F\in \hat{\cal F}\right\}.
\end{equation*}
We denote $\tilde{\cal N}^0$ simply by $\tilde{\cal N}$.

\begin{lemma}\label{lambda-Sigmatilde}
If\, $\Xi_{\tilde {\cal R}}$ is unbounded then $\lambda(\tilde \Sigma)=+\infty$. On the other hand, if $\Xi_{\tilde{\cal R}}$ is bounded, then  
\begin{equation*}
\lambda(\tilde\Sigma^\mu)=
\max\left(\sup_{\s\in \cal K}\alpha(M^\mu),\sup_{\substack{\omega\in \Omega_{\tilde{\cal N}^\mu} \\k\in \mathbb{N}}}\frac{\log(\rho(\Pi_{\omega_k}))}{|\omega_k|}\right)<+\infty,
\end{equation*}
and $\tilde\Sigma^\mu$ 
is ES (respectively, EU) if and only if $\lambda(\tilde\Sigma^\mu)<0$ (respectively, $\lambda(\tilde\Sigma^\mu)>0$).
\end{lemma}

\begin{proof}
The 
proof follows the same lines of that of Lemma~\ref{lambda-Sigmabar}.
\end{proof}

\section{Proof of Theorem~\ref{main theorem-necessary}}\label{sec: necessary}
Let us start this section by stating an approximation result of the exponential of $\Gamma^{\e,\mu}$, whose proof is given in the appendix.

 \begin{lemma}\label{sufficient1}
Let the $D$-Hurwitz assumption holds.
Let $\mu\in\R$ and $\mathcal{T}\subset \R_{\geq 0}$.
Assume that either $\mathcal{T}$ is bounded or $\alpha(M^\mu)<0$ for every $\s\in {\cal K}$. 
Then there exist $C>1$ (independent of $\mu, {\cal T}$) and $K>0$ such that
for every $t\in{\cal T}$, $\s=(\ell,P,\Lambda,R)\in {\cal K}$, and every $\e>0$ small enough 
\begin{itemize}
    \item if\, $t\geq C\e|\log(\e)|$ then  
    \begin{align}\label{estimate11-sufficient}
\left\|(T^{\e})^{-1}e^{t\Gamma^{\e,\mu}}T^{\e}-T^{-1}\left(\begin{smallmatrix} e^{t M^{\mu}}&0\\0& 0\end{smallmatrix}\right)T\right\|\leq K\e; 
\end{align} 
\item if\, $t<C\e|\log(\e)|$ then 
\begin{align}\label{estimate12-sufficient}
\left\|(T^{\e})^{-1}e^{t\Gamma^{\e,\mu}}T^{\e}-T^{-1}\left(\begin{smallmatrix}I_{\ell}&0\\0& e^{\frac{t}{\e}D}\end{smallmatrix}\right)T\right\|\leq Kt,
\end{align}
\end{itemize}
where  $M^{\mu}=M^{\mu}(\s)$ and $\Gamma^{\e,\mu}=\Gamma^{\e,\mu}(\s)$ are defined in~\eqref{Gamma}, and $D=D(\s)$ is given in~\eqref{ABCD}.
\end{lemma}

\subsection{Proof of item~\ref{itemii} of Theorem~\ref{main theorem-necessary}}\label{subsec:itemii}
Let 
$\mu\in \mathbb{R}$ be such that $\mu>-\lambda(\bar\Sigma_{\tau})$, so that $\lambda(\bar\Sigma^{\mu}_{\tau})>0$. 
 From Lemma~\ref{lambda-Sigmabar}, 
there exist either $\sigma\in{\cal K}
 $ such that $\alpha(M^{\mu}(\s))>0$ or $\omega\in\Omega_{{\bar {\cal N}}_{\tau}^{\mu}}$ and $k\in \mathbb{N}$ such that the spectral radius of $\Pi_{\omega_k}$ is greater than one. In the first case, 
by continuity of the spectral abscissa,
we have $\alpha(\Gamma^{\e,\mu})>0$ for sufficiently small $\e>0$.
In the second case, thanks to \eqref{estimate11-sufficient} and the fact that ${\bar {\cal N}}_{\tau}^{\mu}$ is bounded,
there exists $K_\omega>0$ such that for $\e>0$ sufficiently small we have 
\begin{equation*}\label{eq3: comparison lambdas}
\left\|\Pi_{\omega^{\e}_k}-\Pi_{\omega_k}\right\|\leq K_{\omega}\e|\log(\e)|,
\end{equation*}
 where $\omega^\e\in \Omega_{\cal N^{\e,\mu}_{\tau}}$ is the sequence corresponding to $\omega$, in the sense that if the $j$th element of $\omega$ is 
the pair  $(R(\s)T(\s)^{-1}(\begin{smallmatrix}e^{t M^\mu(\s)}&0\\0&0\end{smallmatrix})T(\s),t)$ then 
 the $j$th element of $\omega^\e$ is 
 $(R(\s)T^{\e}(\s)^{-1}e^{t \Gamma^{\e,\mu}}T^{\e}(\s),t)$. 
By consequence, from the continuity of the spectral radius, we have  $\rho\left(\Pi_{\omega^{\e}_k}\right)>1$ for $\e$ small enough. Let $\tilde\omega^\e\in \Omega_{\cal N^{\e,\mu}_{\tau}}$ be $k$-periodic 
such that $\tilde\omega^\e_k$ is given by ${\omega^{\e}_k}$. 
Using Lemma~\ref{lambda-sigmaeps}, it follows that $\lambda(\Sigma^{\e,\mu}_{\cal K,\tau})>0$ for every $\e>0$ small enough. This proves, in particular, the last part of the statement.
Since $\lambda(\bar\Sigma^{\mu}_{\tau})=\mu +\lambda(\bar\Sigma_{\tau})$ and $\lambda(\Sigma^{\e,\mu}_{\cal K, \tau})=\mu+\lambda(\Sigma^{\e}_{\cal K, \tau})$, the inequality stated in item~\ref{itemii} of Theorem~\ref{main theorem-necessary}  is obtained by  letting $\mu\searrow-\lambda(\bar\Sigma_{\tau})$. \\
The rest of the proof follows immediately from Lemma~\ref{lambda-Sigmabar} together with item~\ref{item3} of Lemma~\ref{thm:supsup}.

\subsection{Proof of item~\ref{itemiii} of Theorem~\ref{main theorem-necessary}}\label{subsec:itemiii}
 %\vspace{-0.3cm}
 If $\tau=0$ and both $\Xi_{{\cal R}}$ and $\Xi_{\bar {\cal R}}$ are bounded then $\max\{\lambda(\bar\Sigma),\lambda(\Sigma^{\e}_{\cal K})\}<+\infty$ for every $\e>0$, as it follows from Lemma~\ref{thm:supsup} and Lemmas~\ref{lambda-sigmaeps} and \ref{lambda-Sigmabar}. The remainder of the proof proceeds analogously to that of item~\ref{itemii}.
%\end{proof}

\subsection{Proof of item~\ref{itemv} of Theorem~\ref{main theorem-necessary}}\label{subsec:itemv}

Since $\Xi_{{\cal R}}$ is bounded,  $\lambda(\hat\Sigma)$ is characterized by \eqref{eq:lambda-hatsigma} in Lemma~\ref{lambda-hatsigma}.
Let $\delta>0$. There exist $\bar\omega\in \Omega_{\hat{\cal N}}$ and $j\in \mathbb{N}$ such that 
\begin{equation}\label{eq6: comparison lambdas}
\dfrac{\log(\rho(\Pi_{\bar\omega_j}))}{|\bar\omega_j|}> \sup_{\omega\in \Omega_{\hat{\cal N}},k\in \mathbb{N}}\dfrac{\log(\rho(\Pi_{\omega_k}))}{|\omega_k|}-\delta.
\end{equation}
Let $\bar\omega^\e\in \Omega_{\cal N^\e}$ be the sequence corresponding to $\bar\omega$, in the sense that if the $j$th element of $\bar\omega$ is 
the pair  $(R(\s)T(\s)^{-1}(\begin{smallmatrix}I_{\ell(\s)}&0\\0&e^{s D(\s)}\end{smallmatrix})T(\s),s)$ then 
 the $j$th element of $\bar\omega^\e$ is 
 $(R(\s)T^{\e}(\s)^{-1}e^{s\e \Gamma^{\e,\mu}}T^{\e}(\s),s\e)$. Notice that 
$|\bar\omega^\e_k|=\e|\bar\omega_k|$ for every $k\in \mathbb{N}$.
Thanks to Lemma~\ref{sufficient1} and the continuity of the spectral radius, for $\e$ small enough it holds that 
\begin{equation}\label{eq4: comparison lambdas}
\dfrac{\log(\rho(\Pi_{\bar\omega^\e_j}))}{|\bar\omega_j|}> \dfrac{\log(\rho(\Pi_{\bar\omega_j}))}{|\bar\omega_j|}-\delta.
\end{equation}
Using the fact that $\ds\sup_{\omega\in \Omega_{\cal N^\e},k\in \mathbb{N}}\dfrac{\log(\rho(\Pi_{\omega_k}))}{|\omega_k|}\ge \frac{\log(\rho(\Pi_{\bar\omega^\e_j}))}{\e |\bar\omega_j|}$, we deduce from~\eqref{eq6: comparison lambdas} together with~\eqref{eq4: comparison lambdas} that 
\begin{equation*}
    \e\sup_{\omega\in \Omega_{\cal N^\e},k\in \mathbb{N}}\dfrac{\log(\rho(\Pi_{\omega_k}))}{|\omega_k|}> \sup_{\omega\in \Omega_{\hat{\cal N}},k\in \mathbb{N}}\dfrac{\log(\rho(\Pi_{\omega_k}))}{|\omega_k|}-2\delta. 
\end{equation*}
By arbitrariness of $\delta$, it follows that 
\begin{equation*}\label{eq5: comparison lambdas}
\liminf_{\e\searrow 0}\left\{\e\sup_{\omega\in \Omega_{\cal N^\e},k\in \mathbb{N}}\dfrac{\log(\rho(\Pi_{\omega_k}))}{|\omega_k|}\right\}\ge \sup_{\substack{\omega\in \Omega_{\hat{\cal N}}\\k\in \mathbb{N}}}\dfrac{\log(\rho(\Pi_{\omega_k}))}{|\omega_k|}.
\end{equation*}
By consequence, from Lemmas~\ref{lambda-sigmaeps} and \ref{lambda-hatsigma}  we get 
 \begin{align*} 
 \lambda(\hat\Sigma)&\le \max\left(0, \liminf_{\e\searrow 0}\left\{\e\sup_{\omega\in \Omega_{\cal N^\e},k\in \mathbb{N}}\dfrac{\log(\rho(\Pi_{\omega_k}))}{|\omega_k|}\right\}\right)\\
&\le \max\left(0,\liminf_{\e\searrow 0} \e\lambda(\Sigma^{\e}_{\cal K})\right),
 \end{align*} 
 concluding the proof of the inequality stated in item~\ref{itemv} of Theorem~\ref{main theorem-necessary}. \\
The rest of the proof follows immediately from Lemma~\ref{lambda-hatsigma} together with item~\ref{item3} of Lemma~\ref{thm:supsup}.

\section{Proof of Theorem~\ref{main theorem-sufficient}}\label{sec: sufficient}

We begin this section by stating a converse Lyapunov result which can be obtained by using classical arguments (as, for instance, those of \cite[Lemma~6.2]{cai2007smooth}). For this purpose, we need to consider the quantity $\tilde\lambda(\Delta_{{\cal Z},\tau})$ associated with an  
impulsive linear switched system $\Delta_{{\cal Z},\tau}$, defined as
\[\tilde\lambda(\Delta_{{\cal Z},\tau})=
\limsup_{t\to+\infty}\frac{1}{t}
\sup_{Z\in {\cal S}_{{\cal Z},\tau},\; k\in \Theta^\star(Z),\; t=t_k}\log(\|\Phi_Z(t,0)\|).
\]

\begin{proposition}
\label{converse}
Consider the impulsive linear system  $\Delta_{{\cal Z},\tau}$ with $\tau\ge 0$.
Then
\begin{description}
\item[$(i)$] $\tilde\lambda(\Delta_{{\cal Z},\tau})<0$ if and only if there exist $c>1$, $\gamma>0$, and  $V:\R^{d}\to \R_{\geq 0}$ convex, absolutely homogeneous, and Lipschitz continuous such that, for every $(Z_1,Z_2)\in \mathcal{Z}$, $t\in \mathbb{R}_{\ge\tau}$, $x\in\mathbb{R}^d$, one has 
\begin{align}
|x|\leq V(x) &\leq c|x|,  \label{comp}\\
V(Z_2e^{t Z_1}x)&\leq e^{-\gamma t}V(x). \label{decay}
\end{align}
\item[$(ii)$] $\Delta_{{\cal Z},\tau}$ is ES if and only if 
$\sup_{(Z_1,Z_2)\in\mathcal{Z}}\alpha(Z_1)<0$ and $\tilde\lambda(\Delta_{{\cal Z},\tau})<0$.
\end{description}
\end{proposition}

As a corollary of Proposition~\ref{converse},
we have the following result. 

\begin{corollary}\label{cor-converse}
Let $\tau\ge 0$ and $\mu\in \R$.
Consider one of the following three cases:
\begin{itemize}
\item[(C1)] $\bar\Sigma^{\mu}_\tau$ is ES and ${\cal N}=\bar {\cal N}_\tau^\mu$,
\item[(C2)] $\tilde\Sigma^\mu$ is ES and ${\cal N}=\tilde {\cal N}^\mu$,
\item[(C3)] $\tilde\lambda(\hat\Sigma)<0$ and ${\cal N}=\hat {\cal N}$.
\end{itemize}
Then 
there exist $c>1$, $\gamma>0$, and  $V:\R^{d}\to \R_{\geq 0}$ absolutely homogeneous and Lipschitz continuous 
such that
\eqref{comp} holds true and
\[
V(N x)\leq e^{-\gamma t}V(x), \quad \forall\, x\in \R^{d},
\]
for every $(N,t)\in{\cal N}$.
\end{corollary}
\begin{proof}
In case (C1) observe that, as noticed in Lemma~\ref{lambda-Sigmabar}, the exponential stability of $\bar\Sigma^{\mu}_\tau$ is equivalent to that of $\bar\Delta_{\tau}^{\delta}$, up to replacing $M$ by $M^\mu$ and choosing $\delta<\sup_{\s\in \cal K}\alpha(M^\mu)$.
The conclusion then follows from Proposition~\ref{converse}. The argument for
Case (C2) is similar, with the role of Lemma~\ref{lambda-Sigmabar} played by Lemma~\ref{lambda-Sigmatilde}. As for Case (C3), it is enough to observe that $\tilde\lambda(\hat \Sigma)=\tilde\lambda(\hat \Delta)$, where $\hat\Delta$ is the system introduced in the proof of Lemma~\ref{lambda-hatsigma}, %\eqref{hatDelta}, 
and to apply Proposition~\ref{converse}.
\end{proof}

\subsection{Proof of item~\ref{itemi} of Theorem~\ref{main theorem-sufficient}}\label{subsec:itemi}

In order to prove the inequality stated in item~\ref{itemi} of Theorem~\ref{main theorem-sufficient}, we can assume without loss of generality that $\lambda(\bar\Sigma_{\tau})<+\infty$.
Now, consider $\mu\in \R$ such that $\lambda(\bar\Sigma^{\mu}_{\tau})<0$ (or, equivalently, that 
$\bar\Sigma^{\mu}_{\tau}$ is ES by Lemma~\ref{lambda-Sigmabar})
and let $V$ be the Lyapunov function associated with system $\bar\Sigma^{\mu}_{\tau}$ by Corollary~\ref{cor-converse} in Case~(C1), with corresponding constants $c>1$, $\gamma>0$.
For $\s\in {\cal K}$ and $t\geq\tau$, let $(N^{\e},t)\in {\cal N}^{\e,\mu}_{\tau}$ where $N^{\e}=R(T^{\e})^{-1}e^{t\Gamma^{\e,\mu}}T^{\e}\in {\cal N}^{\e,\mu}_{\tau}$ is the associated evolution.
Since 
$\sup_{\s\in\cal K}\alpha(M^\mu)\leq \lambda(\bar\Sigma^{\mu}_{\tau})$, 
by continuity of the spectral abscissa and compactness of ${\cal K}$,  
it follows 
that $\sup_{\s\in\cal K}\alpha(\Gamma^{\e,\mu})<0$ for every $\e>0$ small enough, which implies
\begin{equation}\label{unif-bound-flows}
\|N^{\e}\|\leq C_0e^{-\eta t},
\end{equation}
for some $C_0>1$ and $\eta>0$ (independent of $\e, \s$). Hence, from~\eqref{comp},
\begin{equation*}\label{suff-dwell-2.0}
V(N^{\e}x)\leq cC_0e^{-\eta t}V(x).
\end{equation*}
If we fix $t^\star>\frac{\log(c C_0)}{\eta}$, we deduce from the inequality above that there exists $\beta>0$ sufficiently small such that, if $t \geq t^\star$ then
\begin{equation}\label{suff-dwell-2.1}
V(N^{\e}x)\leq e^{-\beta t}V(x).
\end{equation}
Furthermore, thanks
to Lemma~\ref{sufficient1}, there exist $K>0$ (independent of $\s, \e$) such that for $\e>0$ small enough we have $\|N^{\e}-N\|\leq K\e$, where $(N,t)$ is the corresponding element in $\bar{\cal N}^\mu_{\tau}$. For $x\in \mathbb{R}^d$, we have   
\begin{align*}
V(N^{\e}x)&\leq L\|N^{\e}-N\||x|+V(Nx)\\
&\leq 
\left(\e KL+e^{-\gamma t}\right)V(x),
\end{align*}
where $L>0$ is such that $V$ is $L$-Lipschitz continuous. 
It follows from the equation above that, for $\e$ sufficiently small and up to reducing $\beta>0$, the inequality~\eqref{suff-dwell-2.1} holds true even if $t\in [\tau,t^\star]$, hence for every $t\in \mathbb{R}_{\ge \tau}$.
By Proposition~\ref{converse}, we deduce then that
$\Sigma^{\e,\mu}_{\tau}$ is ES 
(and hence $\lambda(\Sigma^{\e,\mu}_{\tau})<0$) for $\e$ small enough. 
Since $\lambda(\Sigma^{\e,\mu}_{\cal K,\tau})=\mu+\lambda(\Sigma^{\e}_{\cal K, \tau})$,
and considering the limit as
$\mu\nearrow-\lambda(\bar\Sigma_{\tau})$, 
we deduce that
$\limsup_{\e\searrow0
}\lambda(\Sigma^{\e}_{\cal K,\tau})\leq \lambda(\bar\Sigma_{\tau})$. 
%\end{proof}

\subsection{Proof of item~\ref{itemiv} of Theorem~\ref{main theorem-sufficient}}\label{subsec:itemiv}

First notice that if $\Xi_{\tilde{\cal R}}$ is unbounded, then, by Lemma~\ref{lambda-Sigmatilde}, $\lambda(\tilde\Sigma)=+\infty$ and there is nothing to prove. 
Assume then that $\Xi_{\tilde{\cal R}}$ is bounded.
Let $\mu\in \R$ be such that $\lambda(\tilde\Sigma^{\mu})<0$. 
Let $V$ be the Lyapunov function associated with system $\tilde\Sigma^{\mu}$ by Corollary~\ref{cor-converse} in Case~(C2), 
 with corresponding constants $c>1$ and $\gamma>0$. 
 Let, moreover, $W$ be the Lyapunov function associated with  system $\hat\Sigma$  by Corollary~\ref{cor-converse} in Case~(C3), 
with corresponding constants $c_W>1$ and $\gamma_W>0$.\\
Fix $\e>0$ and $\s\in {\cal S}_{\cal K}$ and associate with them the corresponding sequence $\omega^\e\in \Omega_{{\cal N}^{\e,\mu}
}$. 
Let $(t_i)_{i\in \mathbb{N}}$ be the switching times of $\sigma$, $(\sigma_i)_{i\in \mathbb{N}}$ its switching values, and $(N_i^\e,t_{i+1}-t_i)_{i\in \mathbb{N}}$ the corresponding values of $\omega^\e$. 
Let $C>1$ and $K>0$ be as in Lemma~\ref{sufficient1} with ${\cal T}=\R_{\ge 0}$ (which is possible because 
$\lambda(\tilde\Sigma^{\mu})<0$
implies that $\alpha(M^\mu)<0$ for every $\s\in{\cal K}$, according to Lemma~\ref{lambda-Sigmatilde}).
Fix $\bar s>
\log(c c_W)/\gamma_W>0$, so that, for every $s\ge \bar s$,
\begin{align}
c c_W e^{-\gamma_W s}&
<1. 
\label{eq:lemmahomogeneity}
\end{align}
We say that an interval $[t_i,t_{i+1})$ is of  
\begin{itemize}
\item[]type (a) if $t_{i+1}-t_i\ge C|\e\log(\e)|$, 
\item[]   type (b) if $0<  t_{i+1}-t_i <C|\e\log(\e)|$.
\end{itemize}
For every $i\in \mathbb{N}$, define $N_i$ as follows: if 
$[t_i,t_{i+1})$ is of type (a) (respectively, (b)
), let $(N_i,t_{i+1}-t_i)\in \bar{\cal N}^\mu_0$
(respectively, $(N_i,\frac{t_{i+1}-t_i}{\e})\in\hat{\cal N}$) be the pair corresponding to the same mode $\s_i$ as $N_i^\e$.\\
By Lemma~\ref{sufficient1}
there exists 
$K>0$ such that, for $\e$ small enough, $\|N_i^\e-N_i\|\le K\e$ when $[t_i,t_{i+1})$ is of the type (a) and 
$\|N_i^\e-N_i\|\le K(t_{i+1}-t_i)$ when $[t_i,t_{i+1})$ is of the type (b). 
 \\
Notice that $\hat{\cal N}$ is product bounded, 
as it follows from the inequality $\tilde\lambda(\hat\Sigma)<0$ and from \cite[Lemma 3]{technical-spin-off} 
applied to the system $\hat\Delta$ introduced in the proof of Lemma~\ref{lambda-hatsigma}.
We can then apply Lemma~\ref{eps-estim} given in the Appendix and deduce that there exists 
$\kappa>0$ such that, for $\e$ small enough,  
\[\|N_i^\e\cdots N_j^\e-N_i\cdots N_j \|\le e^{\kappa(t_{i+1}-t_j)}-1\]
 when each of the intervals $[t_j,t_{j+1}),\dots,[t_i,t_{i+1})$  is of the type (b). 
 \\
In order to show that the function $V$ is (exponentially) decreasing along the dynamics whenever $\e$ 
is small enough, we introduce a further splitting of $\mathbb{R}_{\ge 0}$ obtained by suitably regrouping the intervals of type (a) and (b). We say that an interval is of type (B) if it is the maximal union of adjacent intervals of type (b) and its total length is larger than or equal to $\e \bar s$. Let us denote by ${\cal B}^c$ 
the complement in  $\mathbb{R}_{\ge 0}$ of the union of all intervals of type (B).
Then ${\cal B}^c$  is made then of intervals of type (a) and concatenations of intervals of type (b) whose length is smaller than $\e \bar s$.
 Moreover, each bounded connected component of ${\cal B}^c$, except possibly the first one, starts by an interval of type (a), since otherwise the adjacent leftward interval in (B) would not be a maximal union of intervals of type (b). Hence we can split ${\cal B}^c$, except possibly for its first connected component,  in intervals that we call of type (A), namely, 
concatenations of an interval of type (a) followed by  (possibly zero) 
intervals of type (b) of total length smaller than $\e \bar s$.

Fix an interval $[t_{j
},t_{i+1})$ of type (A). 
We distinguish two cases, depending on whether 
$t_{j+1}-t_j$
is larger than a constant $T>0$ to be fixed later. 
Consider first the case where 
$t_{j+1}-t_j\le T$.
Notice that
$(N_i\cdots N_j,t_{j+1}-t_j)$ is an element of $\tilde{\cal N}^\mu$, so that, by hypothesis on $V$, for every $x\in \mathbb{R}^d$, 
\[V(N_i\cdots N_j x)\le e^{-\gamma (t_{j+1}-t_j)}V( x).\]
 Let $L_V>0$ be such that $V$ is $L_V$-Lipschitz continuous.
 Consider also a constant $C_0>0$ such that $N_i^\e$ and every matrix in $\hat{\cal F}$
 has norm smaller than 
 $C_0$.
 Hence,
\begin{align*}
&V(N_i^\e\cdots N_{j}^\e x)\\
&\le L_V\|N_i^\e\cdots N_{j}^\e- N_i\cdots N_{j}\||x|+
V(N_i\cdots N_{j} x)\\
&\le L_V\left(\|N_{i}^\e \cdots N_{j-1}^\e- N_{i}\cdots N_{j-1}\|\| N_{j}^\e\|\right.\\
&+
\|N_{i}\cdots N_{j-1}\| \|N_{j}^\e-N_{j}\| \big)|x|+
V(N_i\cdots N_{j} x)\\
&\le
L_VC_0
(
e^{\kappa
\e \bar s
-1}+K \e)
|x|+
e^{-\gamma(t_{j+1}-t_{j})} V(x).
\end{align*}
By \eqref{comp}, there exists 
  $\tilde K>0$ such that for every $\e>0$ small enough, 
\begin{align*}
V(N_i^\e\cdots N_{j
}^\e x)
&\le \left(\tilde K \e + e^{-\gamma(t_{j+1}-t_{j})}\right) V(x).
\end{align*}  
 Let $\e$ be sufficiently small such that 
$\tilde K\e+e^{-\gamma(t_{j+1}-t_{j})}<e^{-\frac{\gamma}{2}(t_{j+1}-t_{j}+\bar s \e)}$. This choice is possible because for sufficiently small $\e$ we have $s\mapsto f(s)=\tilde K\e+e^{-\gamma s}-e^{-\frac{\gamma}{2} (s+\bar s\e)}<0$ over the interval $[C|\e\log(\e)|,T]$, as it follows form the convexity of $s\mapsto e^{-s}$, which yields that $e^{-s_1}-e^{-s_2}\ge (s_2-s_1)e^{-s_2}$ for $s_1\le s_2$.
  By consequence, for $\e$ sufficiently small we have 
\begin{align}
V(N_i^\e\cdots N_{j
}^\e x)
&\le e^{-\frac{\gamma}{2}(t_{j+1}-t_{j}+\bar s \e)}V(x)\nonumber\\
&\le  e^{-\frac{\gamma}{2}(t_{i+1}-t_{j
})}V(x).\label{without dwell-1}
\end{align}
Consider now the case where $[t_{j},t_{i+1})$ is of type (A) and $t_{i+1}-t_i> T$.
Arguing as in Section~\ref{subsec:itemi}, it follows from 
$\sup_{\s\in\cal K}\alpha(M^\mu)\leq \lambda(\tilde\Sigma^{\mu})<0$
that there exist $\nu,t_*>0$ 
 such that 
$\|N^\e\|\le e^{-\nu\gamma t}$ 
for every 
$t\ge t_*$, where $(N^\e,t)\in{\cal N}^{\e,\mu}
$.
In particular, assuming that $T\ge t_*$, 
\begin{align*}
V(N_i^\e\cdots N_{j}^\e x)\le{}& 
c|N_{i}^\e\cdots N_{j}^\e x|
\\
\le{} &
c\big(\|N_{i}^\e\cdots N_{j-1}^\e- N_{i}\cdots N_{j-1}\|
\\
&
+\|N_{i}\cdots N_{j-1}\|\big)e^{-\nu\gamma (t_{j+1}-t_j)}|x|
\\
\le{}&
c( 
e^{\kappa\e \bar s}-1
+C_0)
e^{-\nu\gamma (t_{j+1}-t_j)}
V(x).
\end{align*}
Up to choosing $T$ large enough and $\e$ small enough, $c( e^{\kappa\e \bar s}-1
+C_0)
e^{-\nu\gamma (t_{j+1}-t_j)}
\le e^{-\frac{\nu\gamma}{2} (
t_{j+1}-t_j+
\bar s\e)}$ for every $t_{j+1}-t_j
> T$, so that
\begin{equation}\label{without dwell-1.5}
V(N_i^\e\cdots N_{j}^\e x)\le e^{-\frac{\nu\gamma}{2} 
(t_{j+1}-t_j)
}
V(x).
\end{equation}
Now, fix an interval $[t_{j},t_{k})$   
of type 
(B).  
Up to splitting  $[t_{j},t_{k})$ 
in sub-intervals, we can assume without loss of generality that 
\[\e \bar s\le t_k-t_j< 2 \e \bar s+C|\e \log(\e)|.\]
Indeed, let $j_0=j$ and $j_1,\dots,j_h$ be the maximal sequence such that   $j_{s+1}$ is the minimal index in $\{j_s+1,\dots,k\}$ such that $t_{j_{s+1}}-t_{j_s}\ge \e \bar s$: then $\e \bar s\le t_{j_{s+1}}-t_{j_s}<  \e \bar s+C|\e \log(\e)|$ for $s=0,\dots,h-1$ and either $j_h=k$ or
$\e \bar s\le t_k-t_{j_{h-1}}< 2 \e \bar s+C|\e \log(\e)|$.\\
Let $L_W>0$ be such that $W$ is $L_W$-Lipschitz continuous. 
Notice that 
\begin{align*}
W(N_{k-1}^\e&\cdots N_{j}^\e x)\\
\le& L_W\|N_{k-1}^\e\cdots N_{j}^\e- N_{k-1}\cdots N_{j}\||x|
\\&
+
W(N_{k-1}\cdots N_{j} x)\\
\le& 
(L_W
(e^{\kappa (t_k-t_j)}-1)
+
e^{-\gamma_W\frac{t_{k}-t_{j}}{\e}}) W(x).
\end{align*}
Recalling that $\bar s$ has been chosen so that \eqref{eq:lemmahomogeneity} holds true,
observe that for  $\e$ sufficiently small 
one has
\begin{align*}
L_W (e^{\kappa t
}
-1)+&e^{-\gamma_W\frac{t}{\e}}
<\frac{e^{-{\gamma_W}{t}}}{c c_W}
\end{align*}
for every $t\in [\e \bar s,2 \e \bar s+C|\e \log(\e)|]$. 
Thus, 
\begin{align*}
W(N_{k-1}^\e\cdots N_{j}^\e x)
&\le 
\frac{
e^{-
{\gamma_W}
{(t_{k}-t_{j})}
}}
{c c_W}
 W(x).
\end{align*}
By consequence, using~\eqref{comp} both for $V$ and $W$ and \eqref{eq:lemmahomogeneity}, it follows that 
\begin{align}
\nonumber
V(&N_{k-1}^\e\cdots N_{j}^\e x)\\
& \le 
c W(N_{k-1}^\e\cdots N_{j}^\e x)\le 
\frac{e^{-\gamma_W(t_{k}-t_{j})}}{c_W}W(x)
 \nonumber\\
&\le e^{-\gamma_W(t_{k}-t_{j})}|x|
\le 
  e^{-\gamma_W(t_{k}-t_{j})}V(x). 
\label{without dwell-2}
\end{align}
Inequalities~\eqref{without dwell-1}, \eqref{without dwell-1.5}, and  \eqref{without dwell-2} imply that $V$ decreases exponentially along the sequence of extremal points of the intervals of type (A) and (B). Using that the flow matrices of $\Sigma^{\e,\mu}$ are uniformly bounded
and that 
$\alpha(\Gamma^{\e,\mu}(\s))<0$ for $\s\in{\cal K}$ and $\e>0$ small, as it follows from $\alpha(M^{\mu}(\s))<0$ (Lemma~\ref{lambda-Sigmatilde}) and $\lim_{\e\searrow 0}\alpha(\Gamma^{\e,\mu}(\s))=\alpha(M^{\mu}(\s))$,
 we deduce that $\Sigma^{\e,\mu}$ is ES.
\\
The desired inequality stated in item~\ref{itemiv} of Theorem~\ref{main theorem-sufficient}  is obtained by taking the limit as $\mu\nearrow -\lambda(\tilde \Sigma)$. This concludes the proof of item~\ref{itemiv} of Theorem~\ref{main theorem-sufficient}.
%\end{proof}

\section{Applications}\label{sec:applications}
 %\vspace{-0.3cm}
\subsection{The complementary case}\label{sec:comp}
 %\vspace{-0.3cm}

In this section, we consider the complementary case, in which system $\Sigma^{\e}_{\tau}$ results from switching between two linear $d$-dimensional systems, the second one obtained by exchanging the slow and fast dynamics of the first system. 
We derive a simple necessary condition for stability when $d \ge 2$, which is also sufficient in the particular case $d=2$. This is formalised in the following proposition.

\begin{proposition}\label{prop: syslimit}
Consider the switched system $\Sigma^{\e}_{\tau}$ defined by the switching under a dwell-time constraint $\tau> 0$ between 
\begin{align*}\label{JJC}
  & 
  \begin{cases}
\dot x = M_{11}(t)x +M_{12}(t)y,\\
   \e\dot y = M_{21}(t)x+M_{22}(t)y,\\
  \end{cases} 
  & 
 \hspace{-0.5cm} \begin{cases}
 \e \dot x =  M_{11}(t)x +M_{12}(t)y, \\
\dot y = M_{21}(t)x+ M_{22}(t)y,
  \end{cases}
\end{align*}
where $x\in \R^\ell$, $y\in \R^{d-\ell}$ with $\ell\in \{ 1,\dots,d-1\}$, and  $M=\left(\begin{smallmatrix}M_{11} & M_{12}\\
M_{21} & M_{22}\end{smallmatrix}\right)\in {\cal S}_{\cal M, \tau}$, ${\cal M}$ being a bounded subset of $\mathbb{R}^{d}$. 
One has the following:
\begin{description}
\item[$(i)$] System~$\Sigma^{\e}_{\tau}$  is EU for every $\tau>0$ and every $\e>0$ small enough if either their exists $M\in {\cal M}$ such that $\max\{\alpha(M_{11}),\alpha(M_{22})\}>0$ or, in the case where the $D$-Hurwitz assumption is satisfied,  
there 
exist $M,N\in {\cal M}$ such that $\rho(M_{11}^{-1}M_{12}N_{22}^{-1}N_{21} )>1$.
\item[$(ii)$] Conversely, in the case when $d=2$ and $\ell=1$, if the $D$-Hurwitz assumption is satisfied (i.e., $M_{11},M_{22}<0$ for every $M\in {\cal M}$) and
${\rho(M_{11}^{-1}M_{12}N_{22}^{-1}N_{21})<1}$ (i.e., 
$
|M_{12}N_{21}|<M_{11}N_{22}
$) for every $M,N\in {\cal M}$, 
then $\Sigma^{\e}_{\tau}$ is ES for every $\tau>0$ and every $\e>0$ small enough.
\end{description}
\end{proposition}

\begin{proof} 
System~$\Sigma^{\e}_{\tau}$ can be equivalently written as system $\Sigma^{\e}_{\cal K,\tau}$ where in this case the compact set ${\cal K}$ is given by 
\begin{align*}
{\cal K}={}&
\{\ell\}\times \{I_{d}\} \times {\cal M} \times \{I_d\}
\cup \{ d-\ell\}\times \{ J_{d}\} \times \{J_d{\cal M}\} \times \{I_d\},
\end{align*}
where $J_{d}=(\begin{smallmatrix}0&I_{d-\ell}\\I_{\ell}&0\end{smallmatrix})$  and
$J_d{\cal M}=\{J_d M\mid M\in {\cal M}\}$.
The first part of item $i)$ is a direct consequence of Proposition~\ref{limit-criterium}. Concerning the second part of point $i)$, 
in this case one can easily verify that $\Xi_{\bar{\cal R}}$, where ${\bar {\cal R}}$ is given by 
\[\bar {\cal R}=\left\{\left(\begin{smallmatrix}
I_{\ell} & 0 \\ 
-M_{22}^{-1}M_{21}& 0
\end{smallmatrix}\right),
\left(\begin{smallmatrix}
0& -M_{11}^{-1}M_{12} 
\\ 0 & I_{d-\ell}
\end{smallmatrix}\right)\mid M\in {\cal M}\right\},
\]
is unbounded, and then thanks to Remark~\ref{lemma:checkbar}, $\bar\Sigma_{\tau}$ is EU.  By Theorem~\ref{main theorem-necessary} it follows that $\Sigma^{\e}_{\tau}$ is EU for every $\tau> 0$ and for every $\e>0$ small enough. \\
Concerning the point $(ii)$, one can easily verify in this case that each matrix $M_k$ is Hurwitz. Given that $\bar\Sigma_{\tau}$ is one-dimensional, it is necessary ES, and the conclusion follows from~Theorem~\ref{main theorem-sufficient}.
\end{proof}
\begin{remark}
Note that in general the condition that the spectral radius of  $M_{11}^{-1}M_{12}N_{22}^{-1}N_{21}$ is smaller than one for every $M,N\in{\cal M}$ is not a sufficient condition for $ \Sigma^{\e}_{\tau}$ to be exponentially stable (see \cite[Example 19]{haidar2024necessary}). 
\end{remark}

\subsection{Numerical example}\label{sec:example}
\begin{figure}[!ht]
\centering
  \includegraphics[scale=0.55]{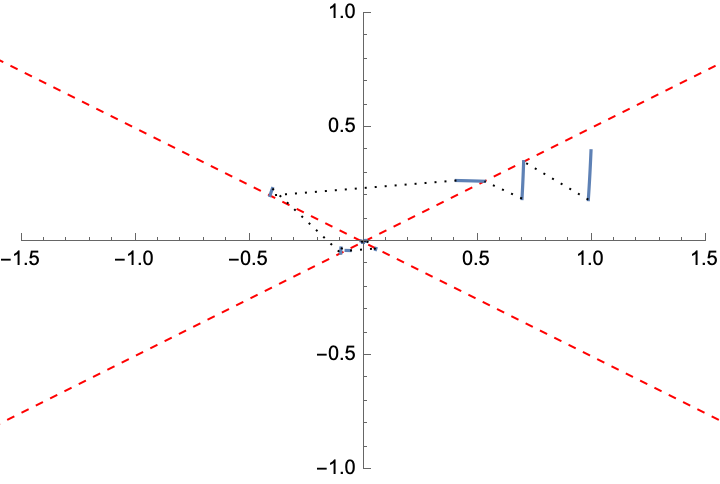}
    \includegraphics[scale=0.55]{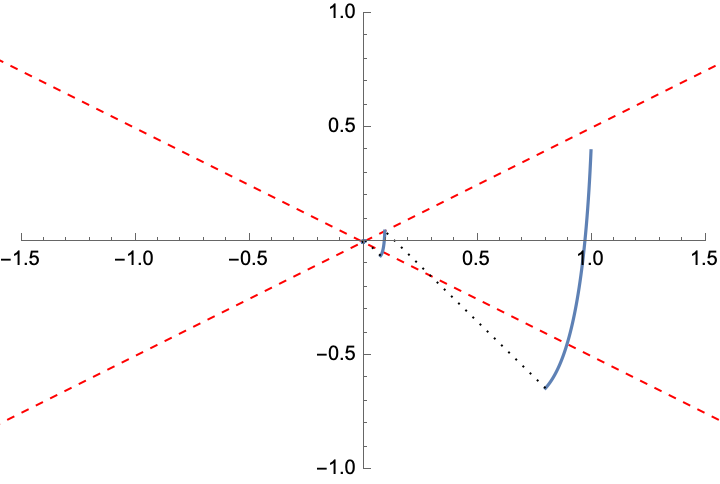}
    \caption{Time evolution
    of system $\Sigma^{\e}_{\cal K}$ of Section~\ref{sec:example} with $r=0.3$, $\e=0.1$, $X(0)=(1,0.4)^T$ corresponding to a
    periodic piecewise-constant switching signal whose period $0.05$ (left) or $0.5$ (right) is divided into three equal intervals corresponding to modes 1, 1, and 2. 
 Dotted lines denote jumps, dashed lines the images of the  jump matrices, i.e., the lines of equation $y=\pm x/2$. 
        }
        
    \label{fig}
\end{figure}

 %\vspace{-0.3cm}
Here, we illustrate through a numerical example the use of the auxiliary system $\tilde\Sigma$ to give a stability criterion for the system $\Sigma^\e_{\cal K}$.\\
For $r\in [0,1]$, consider the switched system $\Sigma^\e_{\cal K}$ in the case where $d=2$ and 
\[{\cal K} = 
\left\{
\begin{aligned}
&\left(1, \,
\left(\begin{smallmatrix}
1 & 0 \\
0 & 1
\end{smallmatrix}\right), \,
\left(\begin{smallmatrix}
-1 & 1 \\
-1 & -1
\end{smallmatrix}\right), \,
\left(\begin{smallmatrix}
2r & 2r \\
r & r
\end{smallmatrix}\right)\right), \\
&\left(1, \,
\left(\begin{smallmatrix}
0 & 1 \\
1 & 0
\end{smallmatrix}\right), \,
\left(\begin{smallmatrix}
1 & -1 \\
-1 & -1
\end{smallmatrix}\right)
, \,
\left(\begin{smallmatrix}
-2r & -2r \\
r & r
\end{smallmatrix}\right)\right)
\end{aligned}
\right\}.
\]
In this case, one can easily verify that
the D-Hurwitz assumption is verified with $D(\sigma)=-1$ for both values of $\sigma$. 
Moreover, $M(\sigma)=-2$ for both values of $\sigma$ and the matrices $J(k),\tilde J(k)$ appearing in systems $\bar \Sigma$ and $\tilde \Sigma$ are always zero, implying that 
such systems are ES. 
The condition of Remark~\ref{rmk:compati}, namely that the points $x\in \R^2$ such that $\left(
\begin{smallmatrix} 
0 & 0\\ 
C(\s) & D(\s)
\end{smallmatrix}
\right)
P(\s)x=0 $ are in the kernel of $R(\sigma)$ for every $\sigma$, can be easily verified.
In addition, it is not difficult to check that  
$\tilde \lambda(\hat\Sigma)<0$ when $r<1/3$. By consequence, thanks to Theorem~\ref{main theorem-sufficient}, we have that $\Sigma^{\e}_{\cal K}$ is ES for every $\e>0$ small enough (see Fig.~\ref{fig}). 

\section{Conclusion}
 %\vspace{-0.3cm}
This paper develops a comprehensive stability analysis for a class of singularly perturbed impulsive linear switched systems characterized by mode-dependent switching between slow and fast dynamics.
Reduced-order single-scale systems are introduced in order to capture the properties resulting from the interaction of the slow and fast dynamics, as  the singular perturbation parameter $\e$ approaches zero. More precisely, the paper 
establishes upper and lower bounds on the maximal Lyapunov exponent of the original system, expressed in terms of the maximal Lyapunov exponents of these auxiliary systems, as $\e$ tends to zero. 
As a consequence,
necessary and sufficient conditions are derived for the exponential stability of the singularly perturbed system for $\e$ small enough. 
Furthermore, a complete characterization of the exponential stability is obtained under a dwell-time constraint on the switching laws.

\section{Appendix: approximations of the flow of singularly perturbed systems 
}
 %\vspace{-0.3cm}
We are interested in this section in approximating the flow of $\Sigma^{\e,\mu}_{\cal K}$
on an interval where the signal $\sigma$ is constant, that is, in approximating $e^{t \Gamma^{\e,\mu}(\s)}$ for some $\s\in{\cal K}$. 
We start our analysis by a useful Gr\"onwall's type result.
\begin{lemma}\label{1xprtoute}
Let $A\in M_n(\R)$ be Hurwitz and $B\in L^\infty([t_0,\infty),M_n(\R))$. 
Then there exist $\alpha,\delta,K,\e_0>0$ depending continuously on $A$  and $\|B\|_{\infty}$ such that, for every 
$\tilde A\in M_n(\R)$ with $\|A-\tilde A\|<\delta$ and every $\e\in (0,\e_0)$
the solution of 
\begin{equation}\label{general-gronwall}
\dot z(t)=\left(\frac{\tilde A}{\e}+ B(t)\right)z(t), \quad \forall\, t\geq t_0, 
\end{equation} 
satisfies the inequalities $|z(t)|\leq Ke^{-\frac{\alpha}{\e}(t-t_0)}|z(t_0)|$
and $|z(t)-e^{\frac{\tilde A}{\e}(t-t_0)}z(t_0) | \leq K \min(\e,t-t_0)|z(t_0)|
$ for every $t\ge t_0$.
\end{lemma}

\begin{proof} 
 By applying the variation of constants formula to~\eqref{general-gronwall}, we obtain
\begin{equation}\label{vc}
z(t)=e^{\frac{\tilde A}{\e}(t-t_0)}z(t_0)+\ds\int_{t_0}^{t}e^{\frac{\tilde A}{\e}(t-s)}B(s)z(s)ds.
\end{equation}
Thanks to the fact that $A$ is Hurwitz, for $\delta$ small enough, there exist $c,\alpha>0$ 
 such that 
\begin{equation}\label{hur}
\|e^{\frac{\tilde A}{\e}(t-s)}\|\leq  ce^{-\frac{2
\alpha}{\e} (t-s)}, \quad \forall\, t\geq s.
\end{equation}
From~\eqref{hur} together with~\eqref{vc} we get that
\begin{equation*}
|z(t)|\leq ce^{-\frac{2
\alpha}{\e}(t-t_0)}|z(t_0)|+c\|B\|_{\infty}\ds\int_{t_0}^{t}e^{-\frac{2
\alpha}{\e}(t-s)}|z(s)|ds,
\end{equation*}
for every $t\geq t_0$, that is, 
\begin{equation*}
\zeta(t)\leq c|z(t_0)|+c\|B\|_{\infty}\ds\int_{t_0}^{t}\zeta(s)ds,\quad \forall\, t\geq t_0,
\end{equation*}
with $\zeta(t)=e^{\frac{2\alpha}{\e}(t-t_0)}|z(t)|$. 
By using Gr\"onwall's inequality, there exists $\e_0$ so that, for every $\e\in(0,\e_0)$, we get that 
\begin{equation}\label{eq:estsed}
|z(t)| 
\leq ce^{-\frac{\alpha}{\e}(t-t_0)}|z(t_0)|,\quad \forall\, t\geq t_0.
\end{equation}
The second inequality in the statement hence follows by 
bounding the integral term in \eqref{vc} using \eqref{hur} and \eqref{eq:estsed}. We first get
\begin{align*}
&|z(t)-e^{\frac{\tilde A}{\e}(t-t_0)}z(t_0)|\leq c^2 \|B\|_\infty (t-t_0)e^{-\frac{\alpha}{\e}(t-t_0)}|z(t_0)|
\\
&\leq c^2 \|B\|_\infty \min\{t-t_0,\e \sup_{s\in\R_{\ge0}} se^{-\alpha s}\}|z(t_0)| ,\quad\forall\;t\ge t_0,
\end{align*}
and, since $s\mapsto s e^{-\alpha s}$ is uniformly bounded on $\R_{\ge0}$, this yields the conclusion.
\end{proof}

We can now state the following. 

\begin{lemma}\label{lemma comparison}
Let the $D$-Hurwitz assumption hold.
Let $\mu\in\R$ and $\mathcal{T}\subset \R_{\geq 0}$.
Assume that either $\mathcal{T}$ is bounded or $\alpha(M^\mu)<0$ for every $\s\in {\cal K}$. 
Then there exists $K>0$ such that 
for 
$(\ell,P,\Lambda,R)\in{\cal K}$,   $t\in \mathcal{T}$, and  $\e>0$ small enough,
\begin{equation}\label{estimate2}
\left\|e^{t\Gamma^{\e,\mu}}-
\begin{pmatrix}
e^{t M^\mu} & 0\\
0& e^{t\frac{D}{\e}}
\end{pmatrix}\right\|\leq K\min(\e,t),
\end{equation}
where  $M^{\mu}$ and $\Gamma^{\e,\mu}$ are defined in~\eqref{Gamma}, and $D=D(\s)$ is given in~\eqref{ABCD}. 
\end{lemma}

\begin{proof}
Let 
$\mu\in\R$  and $(x_0,z_0)\in\R^{\ell}\times\R^{d-\ell}$ be fixed. Consider the trajectory $t\mapsto (x(t),z(t))^T=e^{t \Gamma^{\e,\mu}}(x_0,z_0)^T$. 
As proved in Lemma~\ref{1xprtoute}, there exist $K,\alpha>0$ independent of $(\ell,P,\Lambda,R)\in {\cal K}$ such that  
\begin{equation}\label{eq:zsue}
|z(t)|\leq Ke^{-\frac{\alpha}{\e} t}|z_0|, \qquad \mbox{for }t\ge0,
\end{equation}
 and  
\[|z(t)-e^{\frac t\e D}z_0|\le K\min(\e,t)|z_0|,\qquad \mbox{for }t\ge 0.\]
By a slight abuse of notation, in what follows we still use $K$ to denote possibly larger constants 
independent of $(\ell,P,\Lambda,R)$ and $\e$.
Using that estimate in the dynamics of $x$, we deduce by a simple application of Gr\"onwall's lemma that 
\begin{equation}\label{eq:wiK}
|x(s)|\leq K|(x_0,z_0)|,\qquad \forall\, s\in [0,\sup \mathcal{T}).
\end{equation}
 By applying the variation of constant formula,
we have 
\begin{eqnarray*}
x(t)&=e^{tM^{\mu}}x_0-\e\ds\int_{0}^{t}
e^{(t-s)M^{\mu}}BQ^{\e}x(s)ds\\
&+\ds\int_{0}^{t}e^{(t-s)M^{\mu}}Bz(s)ds.
\end{eqnarray*}
Notice that, by \eqref{eq:zsue} and \eqref{eq:wiK}, we have that
\[
\left|\ds\int_{0}^{t}e^{(t-s)M^{\mu}}Bz(s)ds\right|\leq  K\min(\e,t) |z_0|,
\]
and
\[
\left|\e\ds\int_{0}^{t}
e^{(t-s)M^{\mu}}BQ^{\e}x(s)ds\right|\le K\e t |(x_0,z_0)|.\]

Hence, inequality~\eqref{estimate2} holds. 
\end{proof}

\begin{proof}[Proof of Lemma~\ref{sufficient1}] 
First observe that
 \begin{align*}
T^{\e}&=
\left(\begin{smallmatrix}
I_{\ell}& 0\\
D^{-1}C+\varepsilon Q^{\e} & I_{d-\ell}
\end{smallmatrix}\right)P=T+O(\e), \\
(T^{\e})^{-1}&
=
P^{-1}\left(\begin{smallmatrix}
I_{\ell}& 0\\
-D^{-1}C-\varepsilon Q^{\e} & I_{d-\ell}
\end{smallmatrix}\right) 
=
T^{-1}+O(\e). 
\end{align*}
    Hence,  by the uniform boundedness of $e^{t\Gamma^{\e,\mu}}$ 
    (see equations~\eqref{eq:zsue}-\eqref{eq:wiK}),
\begin{eqnarray}
    (T^{\e})^{-1}e^{t\Gamma^{\e,\mu}}T^{\e}= T^{-1}e^{t\Gamma^{\e,\mu}}T+O(\e).
\label{eq1:eps second order}
  \end{eqnarray} 
 
 By the $D$-Hurwitz assumption, there exist  $c\geq 1$ and $\gamma>0$ depending only on ${\cal K}$ such that 
$\|e^{sD}\|\leq ce^{-\gamma s}$ for all $s\geq 0$. Set $C=\max\{1,1/\gamma\}$ and let us consider  the two  cases   $t\geq C\e|\log(\e)|$ and $t< C\e|\log(\e)|$.\\ 
If $t\geq C\e|\log(\e)|$ then $\|e^{\frac{t}{\e}D}\|\leq ce^{-\gamma \frac{t}{\e}}\leq c\e^{\gamma C}$, and by consequence, thanks to \eqref{estimate2}, one has that 
%\begin{equation*}\label{eq:estimate11}
$\|e^{t\Gamma^{\e,\mu}}-\left(\begin{smallmatrix}e^{tM^{\mu}}&0\\0& 0\end{smallmatrix}\right)\|=O(\e)$.
%\end{equation*}
Inequality~\eqref{estimate11-sufficient} follows from \eqref{eq1:eps second order}.\\

In order to prove the inequality~\eqref{estimate12-sufficient}, we observe, 
as a consequence of the  $D$-Hurwitz assumption, that 
\[\|e^{\frac{t}{\e} D}-I_{d-\ell}\| = O\left(\min\left\{1,\frac{t}{\e}\right\}\right),\] 
from which we deduce that 
\begin{align*}
 &\left\|(T^{\e})^{-1}\left(\begin{smallmatrix}I_{\ell}&0\\0& e^{\frac{t}{\e}D}\end{smallmatrix}\right)T^{\e}-T^{-1}\left(\begin{smallmatrix}I_{\ell}&0\\0& e^{\frac{t}{\e}D}\end{smallmatrix}\right)T\right\|
 \\&=
 \left\|(T^{\e})^{-1}\left(\begin{smallmatrix}0&0\\0& \ e^{\frac{t}{\e}D}\!-I_{d-\ell}\end{smallmatrix}\right)T^{\e}-T^{-1}\left(\begin{smallmatrix}0&0\\0& \ e^{\frac{t}{\e}D}\!-I_{d-\ell}\end{smallmatrix}\right)T\right\|\\
& = O( \min\left\{\e,t\right\}).
\end{align*}
From the estimate above and~\eqref{estimate2} we then have that 
\begin{align*}\label{eq2:eps second order}
 &\left\|(T^{\e})^{-1}e^{t\Gamma^{\e,\mu}}T^{\e}-T^{-1}\left(\begin{smallmatrix}I_{\ell}&0\\0& e^{\frac{t}{\e}D}\end{smallmatrix}\right)T\right\|
 \\&=
\left\| (T^\e)^{-1}\Big(
 e^{t\Gamma^{\e,\mu}}-\left(\begin{smallmatrix}I_{\ell}&0\\0& e^{\frac{t}{\e}D}\end{smallmatrix}\right)
 \Big)T^\e\right\| + O(\min\{\e,t\})\\
 &= \left\| (T^\e)^{-1}
 \left(\begin{smallmatrix}e^{tM^\mu}-I_{\ell}&0\\0& 0\end{smallmatrix}\right)
 T^\e\right\| + O(\min\{\e,t\})\\
 &=O(t),
\end{align*}
and the conclusion follows. 
\end{proof}

\begin{lemma}\label{eps-estim}
Fix an integer $d\ge 1$ and let 
${\cal M}\subset M_{d\times d}(\R)$ 
be product bounded.  
Given $\bar \kappa>0$, there exists a constant $\kappa>0$ such that for every $n\in \mathbb{N}$, 
$t_1, \dots, t_n\in \R_{\ge 0}$, and 
\[
\Psi_1, \dots, \Psi_n\in {\cal M}, 
\quad \Psi'_1, \dots, \Psi'_n\in {\cal M}
\]
satisfying 
\begin{equation}\label{t-estim}
\|\Psi'_j-\Psi_j\|\le \bar \kappa\, t_j, \quad \forall\, 1\le j\le n,
\end{equation} 
we have
\begin{equation}\label{expbound}
\bigl\|\Psi'_1\cdots \Psi'_n - \Psi_1\cdots \Psi_n\bigr\|
   \;\leq\; 
  e^{\kappa(t_1+\cdots+t_n)}-1.
\end{equation}
\end{lemma}

\begin{proof}
Observe that 
\begin{align*}
&\Psi_{1}^{'}\cdots \Psi_{n}^{'}=
\left(\Psi_{1}^{'}-\Psi_{1}+\Psi_{1}\right)\cdots \left(\Psi_{n}^{'}-\Psi_{n}+\Psi_{n}\right)\\
&=\sum_{i\in\{0,1\}^n}
\left(\Psi_{1}^{'}-\Psi_{1}\right)^{i_1}\Psi_1^{1-i_1}\cdots \left(\Psi_{n}^{'}-\Psi_{n}\right)^{i_n}\Psi_n^{1-i_n}.
\end{align*}
For $1\le m\le n$, let $S_m\subset\{0,1\}^n$ denote the set of multi-indices $i$ with exactly $m$ entries equal to $1$.
Then 
\begin{align*}
&\Psi_1^{'}\cdots \Psi_n^{'} -\Psi_1\cdots \Psi_n=\\
&\sum_{m=1}^{n}\sum_{i\in S_m}
\left(\Psi_{1}^{'}-\Psi_{1}\right)^{i_1}\Psi_1^{1-i_1}\cdots \left(\Psi_{n}^{'}-\Psi_{n}\right)^{i_n}\Psi_n^{1-i_n}.
\end{align*}
Observe that each term in the sum is the product of elements of type $\left(\Psi_{j}^{'}-\Psi_{j}\right)\cdots \left(\Psi_{k}^{'}-\Psi_{k}\right)$, for some $1\leq j\leq k\leq n$, and at most $m+1$ elements of type $\Psi_{i}\cdots \Psi_{\ell}$, for some $1\leq i\leq \ell\leq n$. By~\eqref{t-estim}, we have
\begin{align*}
\|\left(\Psi_{j}^{'}-\Psi_{j}\right)\cdots \left(\Psi_{k}^{'}-\Psi_{k}\right)\|\leq (\bar \kappa t_j)\cdots (\bar \kappa t_k).
\end{align*} 
In addition, since ${\cal M}$ is product bounded, there exists $c>0$ such that
\begin{align*}
\|\Psi_{i}\cdots \Psi_{\ell}\|\leq c, \quad \forall\, 1\leq i\leq \ell\leq n. 
\end{align*}
Thus, for each $i\in S_{m}$, 
\begin{align*} 
\left\|\left(\Psi_{1}^{'}-\Psi_{1}\right)^{i_1}\Psi_1^{1-i_1}\cdots \left(\Psi_{n}^{'}-\Psi_{n}\right)^{i_n}\Psi_n^{1-i_n}\right\|\\
\leq  
\bar c^{m+1}\bar\kappa^m\prod_{j=1}^{n}t^{i_j}_j,
\end{align*}
where $\bar c=\max\{1,c\}$. 
Summing over $i\in S_{m}$ gives
\begin{align*} 
\left\|\sum_{i\in S_m}\left(\Psi_{1}^{'}-\Psi_{1}\right)^{i_1}\Psi_1^{1-i_1}\cdots \left(\Psi_{n}^{'}-\Psi_{n}\right)^{i_n}\Psi_n^{1-i_n}\right\|\\
\leq  
\bar c^{m+1}\bar\kappa^m\sum_{i\in S_m}\prod_{j=1}^{n}t^{i_j}_j\le 
(\bar c^{2}\bar\kappa)^m\frac{(t_1+\dots+t_n)^m}{m!}.
\end{align*}
Hence,
\begin{align*} 
\|\Psi_1^{'}\cdots \Psi_n^{'} -&\Psi_1\cdots \Psi_n\|
\leq \sum_{m=1}^{n}\frac{\left(\bar c^2\bar\kappa(t_1+\dots+t_n)\right)^m}{m!}\\
&\leq \sum_{m=0}^{+\infty}\frac{\left(\bar c^2\bar \kappa(t_1+\dots+t_n)\right)^m}{m!}-1\\&
=e^{\bar c^2\bar\kappa(t_1+\dots+t_n)}-1,
\end{align*}
and \eqref{expbound} holds with $\kappa=\bar c^2\bar\kappa$.
\end{proof}

\bibliographystyle{plain}        
\bibliography{biblio}

\end{document}